\documentclass[11pt]{amsart}

\usepackage{geometry}
\usepackage{todonotes}

\usepackage{amsfonts, amssymb, amscd}
\numberwithin{equation}{section}

\usepackage[symbol]{footmisc}

\usepackage{bm}
\usepackage{verbatim}
\usepackage{mathrsfs}
\usepackage{graphicx}
\usepackage{tikz-cd}
\usepackage{subcaption}
\usepackage{listings}
\usepackage{subfiles}
\usepackage[toc,page]{appendix}
\usepackage{mathtools}
\usepackage{comment}
\usepackage{enumerate}
\usepackage{enumitem}
\usepackage[all]{xy}

\usepackage{graphicx}
\graphicspath{{images/}}

\usepackage{appendix}
\usepackage{hyperref}
\hypersetup{
    colorlinks=true,
    citecolor=red,
    linkcolor=blue,
    filecolor=magenta,      
    urlcolor=red,
}
\newcommand{\id}{\mathrm{id}}

\newcommand{\Qq}{\mathbb{Q}}

\newcommand{\Rr}{\mathbb{R}}

\newcommand{\ZZ}{\mathbb{Z}}

\newcommand{\Center}{\operatorname{center}}

\newcommand{\mld}{{\rm{mld}}}

\newcommand{\Supp}{\operatorname{Supp}}

\newcommand{\Bs}{\operatorname{Bs}}
\newcommand{\DG}{\mathcal{DG}}

\newcommand{\mult}{\operatorname{mult}}

\newcommand{\lf}{\lfloor}
\newcommand{\rf}{\rfloor}

\newcommand{\Ii}{\mathcal{I}}

\newtheorem{thm}{Theorem}[section]

\newtheorem{lem}[thm]{Lemma}
\newtheorem{prop}[thm]{Proposition}

\newtheorem{claim}[thm]{Claim}

\theoremstyle{definition}
\newtheorem{defn}[thm]{Definition}
\newtheorem{ques}[thm]{Question}
\theoremstyle{definition}

\newtheorem{ex}[thm]{Example}

\newtheorem{exthm}[thm]{Example-Theorem}

\theoremstyle{definition}

\begin{document}

\title{On the fixed part of pluricanonical systems for surfaces}

\author{Jihao Liu and Lingyao Xie}

\begin{abstract}
We show that $|mK_X|$ defines a birational map and has no fixed part for some bounded positive integer $m$ for any $\frac{1}{2}$-lc surface $X$ such that $K_X$ is big and nef. For every positive integer $n\geq 3$, we construct a sequence of projective surfaces $X_{n,i}$, such that $K_{X_{n,i}}$ is ample, $\mld(X_{n,i})>\frac{1}{n}$ for every $i$, $\lim_{i\rightarrow+\infty}\mld(X_{n,i})=\frac{1}{n}$, and for any positive integer $m$, there exists $i$ such that $|mK_{X_{n,i}}|$ has non-zero fixed part. These results answer the surface case of a question of Xu.
\end{abstract}

\address{Department of Mathematics, Northwestern University, 2033 Sheridan Rd, Evanston, IL 60208}
\email{jliu@northwestern.edu}

\address{Department of Mathematics, The University of Utah, Salt Lake City, UT 84112, USA}
\email{lingyao@math.utah.edu}

\date{\today}

\maketitle

\tableofcontents

\section{Introduction}
We work over the field of complex numbers $\mathbb C$.

\emph{Pluricanonical systems} are central objects in the study of birational geometry. More precisely, given a normal projective variety $X$ such that $K_X$ is effective, we would like to study the behavior of the linear systems $|mK_X|$ for any positive integer $m$.

It is well-known that for any sufficiently divisible $m\gg 0$, the birational map given by $|mK_X|$ is birationally equivalent to the Iitaka fibration of $K_X$.  In 2014, Hacon-M\textsuperscript{c}Kernan-Xu proved that for any lc projective variety $X$ of general type and of fixed dimension, there exists a uniform positive integer $m$ such that $|mK_X|$ defines a birational map \cite[Theorem 1.3]{HMX14} (see also \cite{HM06,Tak06,Tsu99}). In other words, $|mK_X|$ defines a birational morphism $X\backslash\Bs(|mK_X|)\rightarrow\mathbf{P}(|mK_X|)$ for some uniform positive integer $m$, where $\Bs(|mK_X|)$ is the base locus of $|mK_X|$.

It is then natural to ask whether the behavior $|mK_X|$ can be described more accurately. Since we already have a birational morphism $X\backslash\Bs(|mK_X|)\rightarrow\mathbf{P}(|mK_X|)$ for some uniform positive integer $m$, one would like to focus on the asymptotic behavior of $\Bs(|mK_X|)$. As the very first step, we have the following question proposed by Prof. C. Xu to the first author in 2018:
\begin{ques}[Xu]\label{ques: eb codim 2}
Assume that $X$ is a klt projective variety of fixed dimension such that $K_X$ is big and nef. When will we have a uniform positive integer $m$, such that $|mK_X|$ defines a birational map and does not have fixed part?
\end{ques}
Note that it is natural to assume $K_X$ to be nef as we can always run an MMP with scaling and reaches a minimal model for varieties of general type (cf. \cite[Corollary 1.4.2]{BCHM10}).

Question \ref{ques: eb codim 2} naturally arises as a combination of \cite[Theorem 1.3]{HMX14} and the effective base-point-freeness theorem \cite[1.1 Theorem]{Kol93}. Note that when the Cartier index is bounded, $|mK_X|$ not only defines a birational map but is also base-point-free for some uniform positive integer $m$. The interesting cases of Question \ref{ques: eb codim 2} appear when the Cartier index of $K_X$ is unbounded, in which case the uniform base-point-freeness cannot be guaranteed.

Question \ref{ques: eb codim 2} is trivial in dimension $1$ but remained widely open in dimension $\geq 2$. In this paper, we study Question \ref{ques: eb codim 2} when $\dim X=2$. The main theorem of this paper is the following:

\begin{thm}\label{thm: 1/2 lc no fixed part}
There exists a uniform positive integer $m$ satisfying the following. Assume that $X$ is a $\frac{1}{2}$-lc projective surface and $K_X$ is big and nef. Then $|mK_X|$ defines a birational map and does not have fixed part.
\end{thm}

The following Example-Theorem is a complementary statement for Theorem \ref{thm: 1/2 lc no fixed part}, which shows that if the Cartier index of $K_X$ is not bounded and $X$ is not $\frac{1}{2}$-lc, then Theorem \ref{thm: 1/2 lc no fixed part} is not expected to hold.

\begin{exthm}\label{exthm: has fixed part when mld goes to 1/n}
For any integer $n\geq 3$, there exists a sequence of projective surfaces $\{X_i\}_{i=1}^{+\infty}$, such that
\begin{enumerate}
    \item $\mld(X_i)>\frac{1}{n}$ for each $i$ and $\lim_{i\rightarrow+\infty}\mld(X_i)=\frac{1}{n}$,
    \item $K_{X_i}$ is ample, and
    \item if $m_i$ is the minimal positive integer such that $|m_iK_{X_i}|$ defines a birational map and has no fixed part, then $\lim_{i\rightarrow+\infty}m_i=\infty$.
\end{enumerate}
\end{exthm}

Note that the assumptions on $\mld(X)$ in Theorem \ref{thm: 1/2 lc no fixed part} and Example-Theorem \ref{exthm: has fixed part when mld goes to 1/n} are natural assumptions: we are only interested in varieties such that the Cartier index of $K_X$ is not bounded, and if we consider a family of singularities $\{(X\ni x)\}$ such that the index of $K_{X}$ is unbounded, then $\{\mld(X\ni x)\}$ is an infinite set (cf. \cite[Proposition 7.4]{CH21}) and the accumulation points of $\{\mld(X\ni x)\}$ belong to $\{0\}\cup\{\frac{1}{n}|~n\in\mathbb{Z}_{\ge 2}\}$ (cf. \cite[Corollary 3.4]{Ale93}). The $\frac{1}{2}$ accumulation point case is resolved by Theorem \ref{thm: 1/2 lc no fixed part} and the remaining cases are resolved by Example-Theorem \ref{exthm: has fixed part when mld goes to 1/n}.

It is also interesting to ask whether Question \ref{ques: eb codim 2} has a positive answer for canonical or terminal threefolds in dimension $3$, as $1$ is the largest accumulation points of $\mld(X\ni x)$ in dimension $3$ (cf. \cite[Appendix, Theorem]{Sho92}). We will not address this question in this paper, but we will provide a related example (cf. Example-Theorem \ref{ex: dim 3 no fixed curve}).

\medskip

\noindent\textbf{Acknowledgement}.  The authors would like to thank Christopher D. Hacon for useful discussions and encouragements. They would like to thank Chenyang Xu for proposing Question \ref{ques: eb codim 2} to them and sharing useful comments to this question. The authors would like to thank useful discussions with Paolo Cascini, Guodu Chen, Jingjun Han, Junpeng Jiao, Yuchen Liu, Yujie Luo, and Qingyuan Xue. The authors were partially supported by NSF research grants no: DMS-1801851, DMS-1952522 and by a grant from the Simons Foundation; Award Number: 256202.

\section{Preliminaries}
We adopt the standard notation and definitions in \cite{KM98}, and will freely use them.

\begin{defn}[Pairs and singularities]\label{defn: positivity}
    	A pair $(X,B)$ consists of a normal quasi-projective variety $X$ and an $\Rr$-divisor $B\ge0$ such that $K_X+B$ is $\Rr$-Cartier. Moreover, if the coefficients of $B$ are $\leq 1$, then $B$ is called a boundary of $X$.
    	
    Let $E$ be a prime divisor on $X$ and $D$ an $\mathbb R$-divisor on $X$. 	We define $\mult_ED$ to be the multiplicity of $E$ along $D$. Let $\phi:W\to X$
	be any log resolution of $(X,B)$ and let
	$$K_W+B_W:=\phi^{*}(K_X+B).$$
	The \emph{log discrepancy} of a prime divisor $D$ on $W$ with respect to $(X,B)$ is $1-\mult_{D}B_W$ and it is denoted by $a(D,X,B).$  For any positive real number $\epsilon$, we say that $(X,B)$ is lc (resp. klt, $\epsilon$-lc, $\epsilon$-klt) if $a(D,X,B)\ge0$ (resp. $>0$, $\ge \epsilon$, $>\epsilon$) for every log resolution $\phi:W\to X$ as above and every prime divisor $D$ on $W$. We say that $X$ is lc (resp. klt, $\epsilon$-lc, $\epsilon$-klt) if $(X,0)$ is lc (resp. klt, $\epsilon$-lc, $\epsilon$-klt).
	
	A germ $(X\ni x,B)$ consists of a pair $(X,B)$ and a closed point $x\in X$. $(X\ni x,B)$ is called an lc (resp. a klt, an $\epsilon$-lc) germ if $(X,B)$ is lc (resp. klt, $\epsilon$-lc) near $x$. $(X\ni x,B)$ is called $\epsilon$-lc at $x$ if $a(D,X,B)\geq\epsilon$ for any prime divisor $D$ over $X\ni x$ (i.e., $\Center_{X}D=x$).
\end{defn}

\begin{defn}\label{defn: DCC and ACC}
	Let $\Ii$ be a set of real numbers. We say that $\Ii$ satisfies the \emph{descending chain condition} (DCC) if any decreasing sequence $a_1\ge a_2 \ge \cdots \ge a_k \ge\cdots$ in $\Ii$ stabilizes. We say that $\Ii$ satisfies the \emph{ascending chain condition} (ACC) if any increasing sequence in $\Ii$ stabilizes.

\end{defn}

\begin{defn}[Minimal log discrepancies]\label{defn: mld and alct}
Let $(X,B)$ be a pair and $x\in X$ a closed point. The \emph{minimal log discrepancy} of $(X,B)$ is defined as
	$$\mld(X,B):=\inf\{a(E,X,B)\mid E \text{ is an exceptional prime divisor over } X\}.$$
	The \emph{minimal log discrepancy} of $(X\ni x,B)$ is defined as
	$$\mld(X\ni x,B):=\inf\{a(E,X,B)\mid E \text{ is a prime divisor over } X\ni x\}.$$
	If $X$ is $\Qq$-Gorenstein, we define $\mld(X):=\mld(X,0)$. 	If $X$ is $\Qq$-Gorenstein near $x$, we define $\mld(X\ni x):=\mld(X\ni x,0)$.
	For any positive integer $d$, we define
	$$\mld(d):=\{\mld(X\ni x)\mid (X\ni x,0)\text{ is lc, } \dim X=d\}.$$
	\end{defn}

\begin{defn}
Let $X$ be a normal projective variety and $D$ an $\Rr$-divisor on $X$.  We define $$|D|:=\{D'\mid 0\leq D'\sim \lfloor D\rfloor\}.$$
For any $\Rr$-divisor $D$ such that $|D|\not=\emptyset$, the \emph{base locus} of $D$ is
$$\Bs(D):=\cap_{D'\sim D}\Supp D',$$
the \emph{fixed part} of $D$ is the unique $\Rr$-divisor $F\geq 0$, such that
\begin{itemize}
    \item For any $D'\in |D|$, $D'\geq F$, and
    \item $\Bs(|D-F|)$ does not contain any divisor,
\end{itemize}
and the \emph{movable part} of $D$ is $D-F$. We also say that $F$ is the fixed part of $|D|$.

We define $\rho(X)$ to be the Picard number of $X$.
\end{defn}

\begin{defn}
A \emph{surface} is a variety of dimension $2$. A \emph{rational surface} is a projective surface that is birational to $\mathbb P^2$. For ever non-negative integer $k$, the \emph{Hirzebruch surface} $\mathbb F_k$ is $\mathbb P_{\mathbb P^1}(\mathcal{O}_{\mathbb P^1}\oplus\mathcal{O}_{\mathbb P^1}(k))$.
\end{defn}


\begin{defn}\label{defn: crt}
Let $n$ be a non-negative integer, and $C=\cup_{i=1}^nC_i$ a collection of proper curves on a smooth surface $U$. The \emph{determinant} of $C$ is defined as $\det(C):=\det(\{-(C_i\cdot C_j)\}_{1\leq i,j\leq n})$ if $C\not=\emptyset$, and we define $\det(\emptyset)=1$. We define the dual graph $\DG(C)$ of $C$ as follows.
\begin{enumerate}
    \item The vertices $v_i=v_i(C_i)$ of $\DG(C)$ correspond to the curves $C_i$.
    \item For each $i$, $v_i$ is labelled by the integer $e_i:=-(C_i^2)$. $e_i$ is called the \emph{weight} of $v_i$.
    \item For $i\neq j$,the vertices $v_i$ and $v_j$ are connected by $C_i\cdot C_j$ edges.
\end{enumerate}
The \emph{determinant} of $\DG(C)$ is defined as $\det(C)$. For any birational morphism $f: Y\rightarrow X$ between normal surfaces, let $E=\cup_{i=1}^nE_i$ be the reduced exceptional divisor for some non-negative integer $n$. We define $\DG(f):=\DG(E)$. If $f$ is the minimal resolution of $X$ (resp. the minimal resolution of $(X\ni x,0)$ for some closed point $x\in X$), we define $\DG(X):=\DG(f)$ (resp. $\DG(X\ni x):=\DG(f)$).
 \end{defn}

\begin{thm}[{cf. \cite[Theorem 3.2, Corollary 3.4]{Ale93}, \cite{Sho94}}]\label{lem: mld acc dim 2} $\mld(2)$ satisfies the ACC, and the set of accumulation points of $\mld(2)$ is $\{\frac{1}{n}\mid n\geq 2\}\cup\{0\}$.
\end{thm}

\begin{prop}[{cf. \cite[Proposition A.5]{CH21}}]\label{prop: fixed mld imply fixed index}
Let $\Ii_0\subset [0.1]$ be a finite set. Then there exists a positive integer $I$ depending only on $\Ii_0$ satisfying the following. Assume that $(X\ni x,0)$ is an lc surface germ such that $\mld(X\ni x)\in\Ii_0$. Then $IK_X$ is Cartier near $x$.
\end{prop}

\begin{lem}[{\cite[Corollary 2.19]{Ale93}}]\label{lem: upper bound weight epsilon lc}
Let $\epsilon$ be a positive real number and $(X\ni x,0)$ an $\epsilon$-lc surface germ. Then for any vertex $v$ of $\DG(X\ni x)$, the weight of $v$ is $\leq\frac{2}{\epsilon}$.
\end{lem}

\begin{lem}[cf. {\cite[Theorem 3.3]{Ale93},\cite[Lemma A.1]{CH21}}]\label{lem: ale93 lem 3.3}
Let $\epsilon$ be a positive real number. Then there exists a finite set $\mathcal{G}=\mathcal{G}(\epsilon)$ of dual graphs and a finite set $\Ii_0=\Ii_0(\epsilon)$ of positive integers, such that for any $\epsilon$-lc germ $(X\ni x,0)$, one of the following holds:
\begin{enumerate}
    \item $\DG(X\ni x,0)\in\mathcal{G}$.
    \item $\DG(X\ni x,0)$ is of the type as in Figure 1.      \begin{figure}[ht]  
     \begin{tikzpicture}        
         \draw (-1,0) ellipse (1.5 and 0.25);   
         \draw (0.75,0) circle (0.1);
         \node [above] at (0.75,0.2) {\footnotesize$w_1$}; 
         \draw (-0.9,0) ellipse (2 and 0.45);   
         \draw (1.1,0)--(1.4,0);
         \draw (1.5,0) circle (0.1);
         \node [above] at (1.5,0.1) {\footnotesize$2$}; 
         \draw (1.6,0)--(2,0);
         \draw (2.1,0) circle (0.1);
         \node [above] at (2.1,0.1) {\footnotesize$2$};
         \draw (2.2,0)--(2.6,0);
         \draw (2.7,0) circle (0.1);
         \node [above] at (2.7,0.1) {\footnotesize$2$};
         \draw [dashed](2.8,0)--(3.7,0);
         \draw [<-](-2,0.18)--(-3,1);             
         \draw [<-](-0.9,0.45)--(-0.9,1);         
         \node [left] at (-2.9,1) {\footnotesize$q_1$};
         \node [above] at (-0.9,0.9) {\footnotesize$e_1$};
    
         \draw (6.8,0) ellipse (2 and 0.45);         
         \draw (6.9,0) ellipse (1.5 and 0.25);       
         \draw (5.1,0) circle (0.1);  
         \node [above] at (5.1,0.2) {\footnotesize$w_2$}; 
         \draw (4.3,0) circle (0.1);
         \draw (3.8,0) circle (0.1);
         \draw (4.4,0)--(4.8,0);
         \draw (3.9,0)--(4.2,0);
         \node [above] at (4.3,0.1) {\footnotesize$2$}; 
         \node [above] at (3.8,0.1) {\footnotesize$2$}; 
         \node [right] at (8.5,1) {\footnotesize$q_2$};
         \node [above] at (6.8,0.9) {\footnotesize$e_2$};         
         \draw [<-](7.7,0.2)--(8.6,1);    
         \draw [<-](6.8,0.45)--(6.8,1);    
    \end{tikzpicture}
    \caption{}
    \label{fig:1}   
    \end{figure}
    Here $e_1=e_1(X\ni x),q_1=q_1(X\ni x)$ and $e_2=e_2(X\ni x),q_2=q_2(X\ni x)$ are the determinants of the sub-dual graphs, such that $e_1,e_2,q_1,q_2\in\Ii_0$. Moreover, we may assume that
    \begin{enumerate}
        \item either $e_1=w_1=2$ and $q_1=1$, or $w_1>2$; and
        \item either $e_2=w_2=2$ and $q_2=1$, or $w_2>2$.
    \end{enumerate}
    \item $\DG(X\ni x,0)$ is of the type as in Figure 2.
          \begin{figure}[ht]     
    \begin{tikzpicture}
                                                    
         \draw (-1,0) ellipse (1.5 and 0.25);   
         \draw (0.75,0) circle (0.1);
         \draw (-0.9,0) ellipse (2 and 0.45);   
         \draw (1.1,0)--(1.4,0);
         \draw (1.5,0) circle (0.1);
         \node [above] at (1.5,0.1) {\footnotesize$2$}; 
         \draw (1.6,0)--(2,0);
         \draw (2.1,0) circle (0.1);
         \node [above] at (2.1,0.1) {\footnotesize$2$};
         \draw (2.2,0)--(2.6,0);
         \draw (2.7,0) circle (0.1);
         \node [above] at (2.7,0.1) {\footnotesize$2$};
         \draw [dashed](2.8,0)--(3.7,0);
         \draw [<-](-2,0.18)--(-3,1);             
         \draw [<-](-0.9,0.45)--(-0.9,1);         
         \node [left] at (-2.9,1) {\footnotesize$q_1$};
         \node [above] at (-0.9,0.9) {\footnotesize$e_1$};

         \draw (3.9,0) circle (0.1);
         \draw (3.9,0.6) circle (0.1);
         \draw (3.9,-0.6) circle (0.1);
         
         \draw (3.9,0.1)--(3.9,0.5);
         \draw (3.9,-0.1)--(3.9,-0.5);
         
         \node [right] at (4,0) {\footnotesize$2$};
         \node [right] at (4,0.6) {\footnotesize$2$};
         \node [right] at (4,-0.6) {\footnotesize$2$};
    \end{tikzpicture}
     
    \caption{}
    \label{fig:2} 
    \end{figure}  
   Here $e_1=e_1(X\ni x)$ and $q_1=q_1(X\ni x)$ are the determinants of the sub-dual graphs, such that $e_1,q_1\in\Ii_0$.
\end{enumerate}
\end{lem}

\section{Global geometry of smooth surfaces}

\subsection{Some elementary lemmas}

\begin{lem}\label{lem: pe intersection negative self intersection negative}
Let $X$ be a smooth projective surface, $D$ a pseudo-effective $\Rr$-divisor on $X$, and $C$ an irreducible curve on $X$. If $D\cdot C<0$, then $C^2<0$.
\end{lem}
\begin{proof}
Let $D=P+N$ be the Zariski decomposition of $D$ such that $P$ is the positive part and $N$ is the negative part. Since $D\cdot C<0$ and $P$ is nef, $N\cdot C<0$. Since $N\geq 0$, $C\subset\Supp N$ and $C^2<0$.
\end{proof}

\begin{lem}\label{lem: not nef curve is -1}
Let $X$ be a smooth projective surface such that $K_X$ is pseudo-effective. Let $C$ be an irreducible curve on $X$ such that $K_X\cdot C<0$. Then $C^2=K_X\cdot C=-1$. In particular, $C$ is a smooth rational curve.
\end{lem}
\begin{proof}
By Lemma \ref{lem: pe intersection negative self intersection negative}, $C^2<0$. Since $X$ is smooth, $K_X\cdot C\leq-1$ and $C^2\leq -1$. Thus $(K_X+C)\cdot C\leq -2$, which implies that $(K_X+C)\cdot C=-2$,  $C^2=K_X\cdot C=-1$, and  $C$ is a smooth rational curve.
\end{proof}

\begin{lem}\label{lem: smooth psd rational curve negatie self intersection}
Let $X$ be a smooth projective surface such that $K_X$ is pseudo-effective, and $C$ a smooth rational curve on $X$. Then $C^2\leq -1$. 
\end{lem}
\begin{proof}
If not, then $C^2\geq 0$. Since $(K_X+C)\cdot C=-2$, $K_X\cdot C\leq -2<0$. Since $K_X$  is pseudo-effective, $C^2<0$, a contradiction.
\end{proof}

\begin{lem}\label{lem: int 1 contraction still smooth rational curve}
Let $X$ be a smooth projective surface, $C$ an irreducible curve on $X$, $f: Y\rightarrow X$ a blow-up of a closed point, $E$ the exceptional divisor of $f$, and $C_Y$ the strict transform of $C$ on $Y$. If $C_Y\cdot E\leq 1$ and $C_Y$ is a smooth rational curve, then $C$ is a smooth rational curve.
\end{lem}
\begin{proof}
Since $X$ is smooth, $Y$ is smooth. Thus $C_Y\cdot E\in\{0,1\}$. If $C_Y\cdot E=0$, then $f$ is an isomorphism near a neighborhood of $C_Y$ and hence $C$ is a smooth rational curve. If $C_Y\cdot E=1$, then $K_X\cdot C=K_Y\cdot C_Y-1$ and $C^2=C_Y^2+1$, and hence $(K_X+C)\cdot C=(K_Y+C_Y)\cdot C_Y=-2$. Thus $C$ is a smooth rational curve.
\end{proof}

\begin{lem}\label{lem: two -1 curve do not intersect}
Let $X$ be a smooth projective surface such that $K_X$ is pseudo-effective, and $E_1,E_2$ two different smooth rational curves on $X$ such that $E_1^2=E_2^2=-1$. Then $E_1\cdot E_2=0$.
\end{lem}
\begin{proof}
Assume that $E_1\cdot E_2\not=0$, then $E_1\cdot E_2=n\geq 1$ for some positive integer $n$. Let $f: X\rightarrow Y$ be the contraction of $E_1$ and $E_{2,Y}:=f_*E_2$. Then $E_{2,Y}^2=-1+n^2\geq 0$ and $K_Y\cdot E_{2,Y}=-1-n<0$. Since $K_X$ is pseudo-effective, $K_Y$ is pseudo-effective, which contradicts Lemma \ref{lem: pe intersection negative self intersection negative}.
\end{proof}

\begin{lem}\label{lem: 2-1-2 structure does no appear}
Let $X$ be a smooth projective surface such that $K_X$ is pseudo-effective, and $E_1,E_2,E_3$ three different smooth rational curves on $X$. If $E_1^2=E_3^2=-2$ and $E_2^2=-1$, then either $E_1\cdot E_2=0$ or $E_2\cdot E_3=0$.
\end{lem}
\begin{proof}
Assume that $E_1\cdot E_2=n_1>0$ and $E_2\cdot E_3=n_3>0$ for some positive integers $n_1$ and $n_3$. Let $f: X\rightarrow Y$ be the contraction of $E_2$. Then $Y$ is smooth and $K_Y$ is pseudo-effective. Let $E_{1,Y}:=f_*E_1$, and  $E_{3,Y}:=f_*E_3$. Then $E_{1,Y}^2=-2+n_1^2$, $E_{3,Y}^2=-2+n_3^2$, $K_Y\cdot E_{1,Y}=-n_1$, and $K_Y\cdot E_{3,Y}=-n_3$. Thus by Lemma \ref{lem: pe intersection negative self intersection negative}, $n_1=n_3=1$, which implies that $E_1^2=E_3^2=-1$ and $E_1\cdot E_3>0$. By Lemma \ref{lem: int 1 contraction still smooth rational curve}, $E_{1,Y}$ and $E_{3,Y}$ are smooth rational curves, which contradicts Lemma \ref{lem: two -1 curve do not intersect}.
\end{proof}

\begin{lem}\label{lem: k=10-rho}
Let $X$ be a smooth rational surface. Then $K_X^2=10-\rho(X)$.
\end{lem}

\begin{proof}
We may run a $K_X$-MMP 
$f: X:=X_0\xrightarrow{f_1}X_1\xrightarrow{f_2}\dots\xrightarrow{f_n}X_n$
such that either $X_n=\mathbb F_k$ for some non-negative integer $k$ or $X_n=\mathbb P^2$. For any $i\in\{0,1,2,\dots,n-1\}$, we have $K_{X_i}^2=K_{X_{i+1}}^2-1$ and $\rho(X_i)=\rho(X_{i+1})+1$. Thus $K_X^2+\rho(X)=K_{X_n}^2+\rho(X_n)$. If $X_n=\mathbb F_k$ for some non-negative integer $k$, then $K_{X_n}^2+\rho(X_n)=8+2=10$. If $X_n=\mathbb P^2$, then $K_{X_n}^2+\rho(X_n)=9+1=10$. Thus $K_X^2=10-\rho(X)$.
\end{proof}

\subsection{Zariski decomposition}

\begin{lem}\label{lem: zar dec and existed nef divisor}
Let $X$ be a smooth projective surface, and $D$, $\tilde D$ two $\Qq$-divisors on $X$, such that $D\geq\tilde D$ and $\tilde D$ is nef. Let 
$D=P+N$
be the Zariski decomposition of $D$, where $P$ is the positive part and $N$ is the negative part. Then $P\geq\tilde D$.
\end{lem}
\begin{proof}
Assume that $N=\sum_{i=1}^na_iC_i$ and $D-\tilde D=\sum_{i=1}^nb_iC_i+D_0$, where $n$ is a non-negative integer, $C_i$ are distinct irreducible curves, $D_0\geq 0$, and for each $i$, $a_i>0$, $b_i\geq 0$, and $C_i\not\subset\Supp D_0$. Then for every $j\in\{1,2,\dots,n\}$,
\begin{align*}
\sum_{i=1}^na_i(C_i\cdot C_j)&=N\cdot C_j=D\cdot C_j=\tilde D\cdot C_j+(D-\tilde D)\cdot C_j\\
&\geq (D-\tilde D)\cdot C_j=\sum_{i=1}^nb_i(C_i\cdot C_j)+D_0\cdot C_j\geq\sum_{i=1}^nb_i(C_i\cdot C_j),
\end{align*}
which implies that
$\sum_{i=1}^n(a_i-b_i)(C_i\cdot C_j)\geq 0$
for every $j$. Since the intersection matrix $\{(C_i\cdot C_j)\}_{1\leq i,j\leq n}$ is negative definite, $a_i\leq b_i$ for each $i$. Thus $D-\tilde D\geq N$, hence $P\geq\tilde D$.
\end{proof}

\begin{lem}\label{lem: inductive zariski to big and nef}
Let $X$ be a smooth projective surface, $D$ a big Weil divisor on $X$, $\tilde D$ a nef Weil divisor on $X$, and $E$ a Weil divisor on $X$, such that
\begin{itemize}
    \item $D=P+N$ is the Zariski decomposition of $D$, where $P$ is the positive part and $N\geq 0$ is the negative part, 
    \item $E=D-\tilde D\geq 0$, and
    \item $|D|$ defines a birational map.
\end{itemize}
Then there exist a big Weil divisor $D_1$ on $X$ and a Weil divisor $E_1$ on $X$, such that
\begin{enumerate}
\item $D_1=\lfloor P\rfloor$,
    \item $E\geq E_1=D_1-\tilde D\geq 0$,
    \item $|D_1|$ defines a birational map, and
    \item either $N=0$ and $D=P$, or there exists at least one irreducible component $F$ of $\Supp E$ such that $\mult_F(E-E_1)\geq 1$.
\end{enumerate}
\end{lem}
\begin{proof}
We let $D_1:=\lfloor P\rfloor$, then (1) holds. Let $E_1:=D_1-\tilde D$. Since $\tilde D$ is nef and $D\geq\tilde D$, by Lemma \ref{lem: zar dec and existed nef divisor}, $P\geq\tilde D$. Thus $P-\tilde D\geq 0$, and hence $$E_1=D_1-\tilde D=\lfloor P\rfloor-\tilde D=\lfloor P-\tilde D\rfloor\geq 0.$$ Since 
$$E-E_1=D-D_1=P+N-\lfloor P\rfloor=\{P\}+N\geq 0,$$
we deduce (2). Since
$|D_1|=|\lfloor P\rfloor|=|P|\cong |D|,$
$|D_1|$ defines a birational map, hence (3). Finally, if $E-E_1\not=0$, then we are done; otherwise, $E-E_1=0$, hence $\{P\}+N=0$. Thus $N=0$, which implies that $D=P$, hence (4).
\end{proof}

\begin{prop}\label{prop: zariski decomposition give big and nef birational map}
Let $X$ be a smooth projective surface, $D$ a big Weil divisor on $X$, and $\tilde D$ a nef Weil divisor on $X$, such that
\begin{itemize}
    \item $D=P+N$ is the Zariski decomposition of $D$, where $P$ is the positive part and $N\geq 0$ is the negative part, 
    \item $D-\tilde D\geq 0$, and
    \item $|D|$ defines a birational map.
\end{itemize}
Then there exists a Weil divisor $D'$ on $X$, such that
\begin{enumerate}
    \item $D\geq D'\geq\tilde D$, 
    \item $D'$ defines a birational map, and
    \item $D'$ is big and nef.
\end{enumerate}
\end{prop}
\begin{proof}
Let $D_0:=D$, $P_0:=P,N_0:=N$ and $E_0:=D-\tilde D$, and let $r_0$ be the sum of all the coefficients of $E_0$. Then $r_0$ is a non-negative integer.

For any non-negative integer $k$, assume that there exist big Weil divisors $D_1,\dots,D_k$ on $X$, Weil divisors $E_1,\dots,E_k$ on $X$, and non-negative integers $r_1,\dots,r_k$, such that for every $i\in\{0,1,\dots,k\}$,
\begin{itemize}
    \item $D_i=P_i+N_i$ is the Zariski decomposition of $D_i$, where $P_i$ is the positive part and $N_i\geq 0$ is the negative part, 
    \item $E_0\geq E_i=D_i-\tilde D\geq 0$, 
    \item $|D_i|$ defines a birational map, 
    \item $r_k$ is the sum of all the irreducible components of $E_i$ such that $0\leq r_k\leq r_0-k$, and
    \item if $i\geq 1$, then $D_i=\lfloor P_{i-1}\rfloor$.
\end{itemize}
It is clear that these assumptions hold when $k=0$. By Lemma \ref{lem: inductive zariski to big and nef}, there are two cases:

\medskip

\noindent\textbf{Case 1}. $N_k=0$ and $D_k=P_k$. In this case, by our assumptions,
\begin{itemize}
    \item $D_k-\tilde D\geq 0$, hence $D_k\geq\tilde D$,
    \item $E_0\geq D_k-\tilde D$, hence $D\geq D_k$,
    \item $D_k$ is big and defines a birational map, and
    \item $D_k=P_k$ is nef.
\end{itemize}
Thus we may let $D':=D_k$.

\medskip

\noindent\textbf{Case 2}. There exists a big Weil divisor $D_{k+1}$ on $X$, a Weil divisor $E_{k+1}$ on $X$, and a non-negative integer $r_{k+1}$, such that
\begin{itemize}
\item $D_{k+1}=\lfloor P_{k}\rfloor$,
    \item $E_0\geq E_{k+1}=D_{k+1}-\tilde D\geq 0$,
    \item $|D_{k+1}|$ defines a birational map, and
    \item $0\leq r_{k+1}\leq r_k-1$.
\end{itemize}
In this case we may replace $k$ with $k+1$ and apply induction on $k$. Since $0\leq r_k\leq r_0-k$, we have $k\leq r_0$. Thus this process must terminate and we are done.
\end{proof}

\subsection{Effective birationality and existence of special nef \texorpdfstring{$\Qq$}-%
-divisors}

\begin{lem}\label{lem: 192 bpf}
Let $X$ be an lc projective surface such that $K_X$ is big and nef, $f: Y\rightarrow X$ the minimal resolution of $X$, and $E_1,\dots,E_n$ the prime $f$-exceptional divisors. Assume that
$K_Y+\sum_{i=1}^na_iE_i=f^*K_X.$
Then for any positive integer $m$, if there exist integers $r_1,\dots,r_n$, such that
\begin{enumerate}
    \item $0\leq r_i\leq \lfloor ma_i\rfloor$, and
    \item $K_Y+\sum_{i=1}^n\frac{r_i}{m}E_i$ is big and nef,
\end{enumerate}
then $|192mK_X|$ does not have fixed part.
\end{lem}

\begin{proof}
$m(K_Y+\sum_{i=1}^n\frac{r_i}{m}E_i)$ is big and nef and Cartier. By \cite[Theorem 1.1, Remark 1.2]{Fuj09} (see also \cite[1.1 Theorem]{Kol93}), $192m(K_Y+\sum_{i=1}^n\frac{r_i}{m}E_i)$ is base-point-free, which implies that the fixed part of $192m(K_Y+\sum_{i=1}^na_iE_i)$ is supported on $\cup_{i=1}^nE_i$. Thus $|192mK_X|$ does not have fixed part.
\end{proof}

\begin{thm}[c.f. {\cite[Theorem 1.3]{HMX14}}]\label{thm: eb hmx14 surface}
There exists a uniform positive integer $m_1$, such that for any lc surface $X$ such that $K_X$ is big, $|m_1K_X|$ defines a birational map.
\end{thm}

\section{\texorpdfstring{$\frac{1}{3}$}-%
-klt surface surfaces}

\subsection{Classification of \texorpdfstring{$(\frac{1}{3}+\epsilon)$}-%
-lc singularities}
\begin{lem}\label{lem: 1/3 lc must 2 chain and 3}
Let $\epsilon$ be a positive real number. Then there exists a positive integer $n_0=n_0(\epsilon)$ depending only on $\epsilon$ satisfying the following. Assume that $(X\ni x,0)$ a $(\frac{1}{3}+\epsilon)$-lc surface germ. Then
\begin{enumerate}
    \item either $n_0K_X$ is Cartier near $x$, or
    \item $X\ni x$ is a cyclic quotient singularity of type $\frac{1}{2k+1}(1,k)$ for some positive integer $k\geq 10$. In particular, $\DG(X\ni x)$ is the following graph, where there are $k-1$ `` $2$" in the graph.
    	
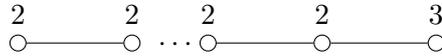
\begin{center}
		\begin{tikzpicture}[scale=1, baseline=(current bounding box.center)]
		\node[circle,draw,inner sep=0pt,minimum size=6pt,label={$2$}] (a) at (0,0){};
		\node[circle,draw,inner sep=0pt,minimum size=6pt,label={$2$}] (b) at (1.5,0){};
		\node[circle,draw,inner sep=0pt,minimum size=6pt,label={$2$},label={[shift={(-0.4,-0.35)}]$\cdots$}] (c) at (2.5,0){};
		\node[circle,draw,inner sep=0pt,minimum size=6pt,label={$2$}] (d) at (4,0){};
		\node[circle,draw,inner sep=0pt,minimum size=6pt,label={$3$}] (e) at (5.5,0){};
		
		\draw (a) edge (b);
		\draw (c) edge (d);
		\draw (d) edge (e);
		\end{tikzpicture}
	\end{center}
\end{enumerate}
\end{lem}
\begin{proof}
Assume that the lemma does not hold. Then there exists a sequence of $(\frac{1}{3}+\epsilon)$-lc surface germs $(X_i\ni x_i,0)$, and a strictly increasing sequence of positive integer $n_i$, such that
\begin{itemize}
    \item $nK_{X_i}$ is not Cartier near $x_i$ for any positive integer $n\leq n_i$, and
    \item $X_i\ni x_i$ is not a cyclic quotient singularity of type $\frac{1}{2k+1}(1,k)$ for any $i$ and any positive integer $k$.
\end{itemize}
We consider the set $\mathcal{A}:=\{\mld(X_i\ni x_i)\}_{i=1}^{+\infty}$. Since $\mld(X_i\ni x_i)\geq\frac{1}{3}+\epsilon$, by Theorem \ref{lem: mld acc dim 2}, the only possible accumulation point of $\mathcal{A}$ is $\frac{1}{2}$. If $\mathcal{A}$ is a finite set, it contradicts Proposition \ref{prop: fixed mld imply fixed index}. Thus possibly passing to a subsequence and replacing $\mathcal{A}$, we may assume that $\mld(X_i\ni x_i)$ is strictly decreasing and $\lim_{i\rightarrow+\infty}\mld(X_i\ni x_i)=\frac{1}{2}$. By Lemma \ref{lem: ale93 lem 3.3}, possibly passing to a subsequence again, we may assume that
\begin{itemize}
    \item either $(X_i\ni x_i)$ satisfies (2) of Lemma \ref{lem: ale93 lem 3.3} for each $i$, and $e_1=e_1(X_i\ni x_i),q_1=q_1(X_i\ni x_i),$ $e_2=e_2(X_i\ni x_i),q_2=q_2(X_i\ni x_i)$ are constants for each $i$, or
    \item $(X_i\ni x_i)$ satisfies (3) of Lemma \ref{lem: ale93 lem 3.3} for each $i$, and $e_1=e_1(X_i\ni x_i),q_1=q_1(X_i\ni x_i)$ are constants for each $i$.
\end{itemize}
We remark that $e_1,e_2,q_1,q_2$ are the numbers defined in Lemma \ref{lem: ale93 lem 3.3}.

If $(X_i\ni x_i)$ satisfies (3) of Lemma \ref{lem: ale93 lem 3.3} for each $i$, and $e_1=e_1(X_i\ni x_i),q_1=q_1(X_i\ni x_i)$ are constants for each $i$, then it contradicts \cite[Lemma A.1(3)]{CH21}. 

If $(X_i\ni x_i)$ satisfies (2) of Lemma \ref{lem: ale93 lem 3.3} for each $i$, and $e_1=e_1(X_i\ni x_i),q_1=q_1(X_i\ni x_i),$ $e_2=e_2(X_i\ni x_i),q_2=q_2(X_i\ni x_i)$ are constants for each $i$, then by \cite[Lemma A.1(3)]{CH21}, $e_1-q_1\leq 2$ and $e_2-q_2\leq 2$. We get a contradiction by enumerating possibilities as follows:

\medskip

\noindent\textbf{Case 1}. $q_1=1$. Then $e_1=2$ or $3$.

\medskip

\noindent\textbf{Case 1.1} $e_1=2$. 

\medskip

\noindent\textbf{Case 1.1.1} $q_2=1$. Then $e_2=2$ or $3$.

\medskip

\noindent\textbf{Case 1.1.1.1} $e_2=2$. In this case, all the weights in $\DG(X_i\ni x_i)$ are $2$. Thus $\mld(X_i\ni x_i)=1$ for every $i$, a contradiction.

\medskip

\noindent\textbf{Case 1.1.1.2} $e_2=3$. In this case, for each $i$, $X_i\ni x_i$ is a cyclic quotient singularity of type $\frac{1}{2k_i+1}(1,k_i)$ for some positive integer $k_i$, and $k_i\rightarrow+\infty$ when $i\rightarrow+\infty$, a contradiction. 

\medskip

\noindent\textbf{Case 1.1.2} $q_2\geq 2$. In this case, there exist an integer $w_2\geq 3$ and a non-negative integer $d_2<q_2$, such that
$e_2=w_2q_2-d_2.$
Thus 
$$2\geq e_2-q_2=(w_2-1)q_2-d_2\geq (w_2-2)q_2+1\geq q_2+1\geq 3,$$
a contradiction.

\medskip

\noindent\textbf{Case 1.2} $e_1=3$.

\medskip

\noindent\textbf{Case 1.2.1} $q_2=1$. Then $e_2=2$ or $3$.

\medskip

\noindent\textbf{Case 1.2.1.1} $e_2=2$. In this case, for each $i$, $X_i\ni x_i$ is a cyclic quotient singularity of type $\frac{1}{2k_i+1}(1,k_i)$ for some positive integer $k_i$, and $k_i\rightarrow+\infty$ when $i\rightarrow+\infty$, a contradiction. 

\medskip

\noindent\textbf{Case 1.2.1.2} $e_2=3$. In this case, $X_i\ni x_i$ is a cyclic quotient singularity of type $\frac{1}{4k_i+8}(1,2k_i+3)$ for some non-negative integer $k_i$, and hence $\mld(X_i\ni x_i)=\frac{1}{2}$, a contradiction.

\medskip

\noindent\textbf{Case 1.2.2} $q_2\geq 2$. In this case, there exists an integer $w_2\geq 3$ and a non-negative integer $d_2<q_2$, such that
$e_2=w_2q_2-d_2.$
Thus 
$$2\geq e_2-q_2=(w_2-1)q_2-d_2\geq (w_2-2)q_2+1\geq q_2+1\geq 3,$$
a contradiction.

\medskip

\noindent\textbf{Case 2}. $q_1\geq 2$. In this case, there exists an integer $w_1\geq 3$ and a non-negative integer $d_1<q_1$, such that
$e_1=w_1q_1-d_1.$
Thus 
$$2\geq e_1-q_1=(w_1-1)q_1-d_1\geq (w_1-2)q_1+1\geq q_1+1\geq 3,$$
a contradiction.
\end{proof}

\subsection{Intersection numbers}

\begin{lem}\label{lem: not rational is pe}
Let $X$ be a projective klt surface such that $K_X$ is nef and $f: Y\rightarrow X$ the minimal resolution of $X$. If $X$ is not rational, then $K_Y$ is pseudo-effective.
\end{lem}
\begin{proof}
If $X$ is not rational, $Y$ is not rational.  If $K_Y$ is not pseudo-effective, then there exists a birational morphism $g: Y\rightarrow W$ to a smooth projective surface $W$ and a $\mathbb P^1$-fibration $h: W\rightarrow R$. Since $Y$ is not a rational surface, $g(R)\geq 0$. Thus for any exceptional curve $F$ of $f$, $f$ does not dominant $R$. Pick a general $g$-vertical curve $\Sigma$ and let $\Sigma_Y$, $\Sigma_X$ be the strict transforms of $\Sigma$ on $Y$ and $X$ respectively. Then
$$0\leq K_X\cdot\Sigma_X=K_Y\cdot\Sigma_Y=K_W\cdot\Sigma=-2,$$
a contradiction.
\end{proof}

\begin{lem}\label{lem: cy only intersect -3 curve}
Let $X$ be a $\frac{1}{3}$-klt surface such that $K_X$ is big and nef, $C$ an irreducible curve on $X$, $x\in C$ a closed point, $f: Y\rightarrow X$ the minimal resolution of $X$, and $C_Y$ the strict transform of $C$ on $Y$. Assume that
\begin{itemize}
\item $X$ is not a rational surface,
    \item $K_Y\cdot C_Y<0$,
    \item $X\ni x$ is a cyclic quotient singularity of type $\frac{1}{2k+1}(1,k)$ for some integer $k\geq 5$, and
    \item $E_1,\dots,E_k$ are prime $f$-exceptional divisors over $X\ni x$, such that
    \begin{enumerate}
        \item $E_i^2=-2$ when $1\leq i\leq k-1$,
        \item $E_k^2=-3$, and
        \item $E_i\cdot E_j\not=0$ if and only if $|i-j|\leq 1$.
    \end{enumerate}
\end{itemize}
Then
\begin{enumerate}
    \item $C_Y\cdot E_i=0$ when $1\leq i\leq k-1$, and
    \item $C_Y\cdot E_k=1$.
\end{enumerate}
\end{lem}

\begin{proof}
By Lemma \ref{lem: not rational is pe}, $K_Y$ is pseudo-effective. By Lemma \ref{lem: not nef curve is -1}, $K_Y\cdot C_Y=-1$ and $C_Y^2=-1$. Moreover, each $E_i$ is a smooth rational curve. We may let $g: Y\rightarrow W$ be the contraction of $C_Y$ and $E_{i,W}:=g_*E_i$ for each $i$. Then $W$ is smooth and $K_W$ is pseudo-effective.

\begin{claim}\label{claim: cy dot ej leq 1}
 $C_Y\cdot E_j\leq 1$ for every $j\in\{1,2,\dots,k\}$.
\end{claim}
\begin{proof}[Proof of Claim \ref{claim: cy dot ej leq 1}]
Suppose not, then there exists an integer $n\geq 2$ an integer $j\in\{1,2,\dots,k\}$, such that $C_Y\cdot E_j=n$. We have
$$E_{j,W}^2=E_j^2+n^2\geq -3+4\geq 1$$
and
$$K_W\cdot E_{j,W}=K_Y\cdot E_j-n\leq -1,$$
which contradicts Lemma \ref{lem: pe intersection negative self intersection negative} as $K_W$ is pseudo-effective. 
\end{proof}

\begin{claim}\label{claim: eiw still rational curve}
$E_{i,W}$ are smooth rational curves for every $i$.
\end{claim}
\begin{proof}[Proof of Claim \ref{claim: eiw still rational curve}]
It immediately follows from Lemma \ref{lem: int 1 contraction still smooth rational curve} and Claim \ref{claim: cy dot ej leq 1}.
\end{proof}

\begin{claim}\label{claim: cy does not intersect middle ej}
$C_Y\cdot E_j=0$ for every $j\in\{2,3,\dots,k-2\}$.  
\end{claim}
\begin{proof}[Proof of Claim \ref{claim: cy does not intersect middle ej}]
Suppose that the claim does not hold. Then by Claim \ref{claim: cy dot ej leq 1}, there exists $j\in\{2,3,\dots,k-2\}$ such that $C_Y\cdot E_j=1$. There are three cases:

\medskip

\noindent\textbf{Case 1}. $C_Y\cdot E_{j+1}=1$. In this case, $E_{j,W}^2=E_{j+1,W}^2=-1$ and $E_{j,W}\cdot E_{j+1,W}=2$. By Claim \ref{claim: eiw still rational curve}, $E_{j,W}^2$ and $E_{j+1,W}$ are smooth rational curves. Since $K_W$ is pseudo-effective, it contradicts Lemma \ref{lem: two -1 curve do not intersect}.

\medskip

\noindent\textbf{Case 2}. $C_Y\cdot E_{j-1}=1$. In this case, $E_{j,W}^2=E_{j-1,W}^2=-1$ and $E_{j,W}\cdot E_{j-1,W}=2$. By Claim \ref{claim: eiw still rational curve}, $E_{j,W}^2$ and $E_{j-1,W}$ are smooth rational curves. Since $K_W$ is pseudo-effective, it contradicts Lemma \ref{lem: two -1 curve do not intersect}.

\medskip

\noindent\textbf{Case 3}. $C_Y\cdot E_{j-1}=C_Y\cdot E_{j+1}=0$. In this case, $E_{j,W}^2=-1$, $E_{j-1,W}^2=E_{j+1,W}^2=-2$, $E_{j,W}\cdot E_{j-1,W}=E_{j,W}\cdot E_{j+1,W}=1$, which contradicts Lemma \ref{lem: 2-1-2 structure does no appear}.
\end{proof}

\begin{claim}\label{claim: cy dot ek-1=0}
$C_Y\cdot E_{k-1}=0$.
\end{claim}
\begin{proof}[Proof of Claim \ref{claim: cy dot ek-1=0}]
Suppose that the claim does not hold. Then by Claim \ref{claim: cy dot ej leq 1}, $C_Y\cdot E_{k-1}=1$. By Claim \ref{claim: cy does not intersect middle ej}, $C_Y\cdot E_j=0$ for every $j\in\{2,3,\dots,k-2\}$. There are two cases:

\medskip

\noindent\textbf{Case 1}. $C_Y\cdot E_k=1$. In this case, $E_{k-1,W}^2=-1$, $E_{k,W}^2=-2$, $E_{k-2,W}^2=-2$, $E_{k-1,W}\cdot E_{k,W}=2$, and $E_{k-1,W}\cdot E_{k-2,W}=1$. This contradicts Lemma \ref{lem: 2-1-2 structure does no appear}. 

\medskip

\noindent\textbf{Case 2}. $C_Y\cdot E_k=0$. In this case, $E_{k-1,W}^2=-1$, $E_{k,W}^2=-3$, $E_{k-2,W}^2=E_{k-3,W}^2=-2$, and for every $i,j\in\{k-3,k-2,k-1,k\}$, $E_i\cdot E_j=1$ if $|i-j|=1$ and $E_i\cdot E_j=0$ if $|i-j|\geq 2$.

Let $h: W\rightarrow Z$ be the contraction of $E_{k-1,W}$ and $E_{i,Z}:=h_*E_{i,W}$ for any $i\not=k-1$. Then $Z$ is smooth and $K_Z$ is pseudo-effective. By Lemma \ref{lem: int 1 contraction still smooth rational curve}, $E_{k-3,Z},E_{k-2,Z}$ and $E_{k,Z}$ are smooth rational curves. Moreover, $E_{k-3,Z}^2=E_{k,Z}^2=-2$, $E_{k-2,Z}^2=-1$, and $E_{k-3,Z}\cdot E_{k-2,Z}=E_{k-2,Z}\cdot E_{k,Z}=1$. This contradicts Lemma \ref{lem: 2-1-2 structure does no appear}.
\end{proof}

\begin{claim}\label{claim: cy cdot e1=0}
$C_Y\cdot E_1=0$.
\end{claim}
\begin{proof}[Proof of Claim \ref{claim: cy cdot e1=0}]
Suppose that the claim does not hold. Then by Claim \ref{claim: cy dot ej leq 1}, $C_Y\cdot E_1=1$. By Claim \ref{claim: cy does not intersect middle ej} and Claim \ref{claim: cy dot ek-1=0}, $C_Y\cdot E_j=0$ for every $j\in\{2,\dots,k-1\}$. By Claim \ref{claim: cy dot ej leq 1}, there are two cases:

\medskip

\noindent\textbf{Case 1}. $C_Y\cdot E_k=1$. In this case, $E_{1,W}^2=1$, $E_{2,W}^2=E_{k,W}^2=-2$, $E_{1,W}\cdot E_{2,W}=E_{1,W}\cdot E_{k,W}=1$, which contradicts Lemma \ref{lem: 2-1-2 structure does no appear}.

\medskip

\noindent\textbf{Case 2}. $C_Y\cdot E_k=0$. The are two sub-cases:

\medskip

\noindent\textbf{Case 2.1}. For any closed point $y\in C$ such that $y\not=x$, $X$ is smooth near $y$. In this case, let $a:=1-a(E_1,X,0)=\frac{1}{2k+1}$. Since $K_X$ is big and nef, 
$$0\leq K_X\cdot C=f^*K_X\cdot C_Y=(K_Y+(1-a)E_1)\cdot C_Y=-1+(1-a)=-a<0,$$
a contradiction.

\medskip

\noindent\textbf{Case 2.2}. There exists a closed point $y\in C$ such that $y\not=x$ and $X$ is not smooth near $y$. Then there exists a prime divisor $F$ on $Y$ that is over $X\ni y$, such that $C_Y\cap F\not=\emptyset$. Moreover, $F$ is a smooth rational curve. Since $X$ is $\frac{1}{3}$-klt, by Lemma \ref{lem: upper bound weight epsilon lc}, $F^2\geq -5$. Let $F_W:=g_*F$. 

We have $F_W\cdot E_{i,W}=0$ for every $i\not=1$, $E_{1,W}^2=-1$, $E_{2,W}^2=E_{3,W}^2=E_{4,W}^2=-2$, and for every $i,j\in\{1,2,3,4\}$, $E_{i,W}\cdot E_{j,W}=1$ when $|i-j|=1$ and $E_{i,W}\cdot E_{j,W}=0$ when $|i-j|\geq 2$.

There are two sub-cases:

\medskip

\noindent\textbf{Case 2.2.1}. $C_Y\cdot F=1$. In this case, by Lemma \ref{lem: int 1 contraction still smooth rational curve}, $F_W$ is a smooth rational curve. Moreover, $F_W^2\geq -4$ and $F_W\cdot E_{1,W}=1$,

Let $h: W\rightarrow Z$ be the contraction of $E_{1,W}$, $E_{i,Z}:=h_*E_{i,W}$ for each $i\not=1$, and $F_Z:=h_*F_W$. Then $Z$ is smooth and $K_Z$ is pseudo-effective. By Lemma \ref{lem: int 1 contraction still smooth rational curve}, $E_{2,Z},E_{3,Z},E_{4,Z}$ and $F_Z$ are smooth rational curves. Moreover, $E_{2,Z}^2=-1$, $E_{3,Z}^2=E_{4,Z}^2=-2$, $E_{2,Z}\cdot E_{3,Z}=E_{3,Z}\cdot E_{4,Z}=F_Z\cdot E_{2,Z}=1$, $F_{2,Z}\cdot E_{3,Z}=F_{2,Z}\cdot E_{4,Z}=E_{2,Z}\cdot E_{4,Z}=0$, and $F_Z^2\geq -3$. 

Let $p: Z\rightarrow T$ be the contraction of $E_{2,Z}$, $E_{i,T}:=p_*E_{i,Z}$ for each $i\not=1,2$, and $F_T:=p_*F_Z$. Then $T$ is smooth and $K_T$ is pseudo-effective. By Lemma \ref{lem: int 1 contraction still smooth rational curve}, $E_{3,T},E_{4,T}$ and $F_T$ are smooth rational curves.  Moreover, $E_{3,T}^2=-1$, $E_{4,T}^2=-2$, $F_T^2\geq -2$, and $E_{3,T}\cdot E_{4,T}=F_T\cdot E_{3,T}=1$. 

By Lemma \ref{lem: smooth psd rational curve negatie self intersection}, $F_T^2\in\{-1,-2\}$. By Lemma \ref{lem: two -1 curve do not intersect}, $F_T^2=-2$. But this contradicts Lemma \ref{lem: 2-1-2 structure does no appear}.

\medskip

\noindent\textbf{Case 2.2.2}. $C_Y\cdot F\geq 2$. In this case, we let $b:=F^2$ and $c:=C_Y\cdot F$. Then $F_W^2=b+c^2$, $K_W\cdot F_W=K_Y\cdot F-c=-2-b-c$, and $F_W\cdot E_{1,W}=c$. 

Let $h: W\rightarrow Z$ be the contraction of $E_{1,W}$ and $F_Z:=h_*F_W$. Then $Z$ is smooth and $K_Z$ is pseudo-effective. Moreover, $F_Z^2=F_W^2+c^2=b+2c^2$, and $K_Z\cdot F_Z=K_W\cdot F_W-c=-2-b-2c$. Since $b\geq -5$ and $c\geq 2$, $F_Z^2\geq 3>0$ and $K_Z\cdot F_Z\leq -1<0$, which contradicts Lemma \ref{lem: pe intersection negative self intersection negative}.
\end{proof}

\noindent\textit{Proof of Lemma \ref{lem: cy only intersect -3 curve} continued}. By Claim \ref{claim: cy does not intersect middle ej}, Claim \ref{claim: cy dot ek-1=0}, and Claim \ref{claim: cy cdot e1=0}, we get (1). Since $x\in C$, $C_Y$ intersects $\cup_{i=1}^kC_i$, which implies that $C_Y$ intersects $E_k$. Thus $C_Y\cdot E_k\geq 1$. (2) follows from Claim \ref{claim: cy dot ej leq 1}.
\end{proof}

\begin{lem}\label{lem: 1/2 lc no k 2k+1 k >10}
Let $X$ be a rational $\frac{2}{5}$-klt surface such that $K_X$ is big and nef and $k\geq 10$ an integer. Then $X$ does not contain a cyclic quotient singularity of type
$\frac{1}{2k+1}(1,k).$
\end{lem}
\begin{proof}
Assume not. Then there exists a closed point $x\in X$ such that $x$ is a cyclic quotient singularity of type $\frac{1}{2k+1}(1,k).$
Let $f: Y\rightarrow X$ be the minimal resolution of $X$ and write
$$K_Y+\sum_{i=1}^k\frac{i}{2k+1}E_i+\sum_{i=1}^sb_iF_i+\sum_{i=1}^tc_iG_i=f^*K_X$$
where $E_1,\dots,E_k,F_1,\dots,F_s,G_1,\dots,G_t$ are the prime $f$-exceptional divisors, where
\begin{itemize}
    \item $E_1,\dots,E_k$ are the prime $f$-exceptional divisors over $X\ni x$,
    \item for every $i\in\{1,2,\dots,s\}$, $\Center_XF_i=x_i$ for some closed point $x_i\in X$, such that $x_i\not=x$ and $\mld(X\ni x_i)\geq\frac{1}{2}$, and
    \item for every $i\in\{1,2,\dots,t\}$, $\Center_XG_i=y_i$ for some closed point $y_i\in X$, such that $y_i\not=x$ and $\mld(X\ni y_i)<\frac{1}{2}$.
\end{itemize}
Since $X$ is $\frac{2}{5}$-klt, by Lemma \ref{lem: upper bound weight epsilon lc}, each $-F_i^2\leq 4$ and $-G_i^2\leq 4$. Moreover, $y_i$ is a cyclic quotient singularity of type $\frac{1}{7}(1,2)$ for every $i$. Since the dual graph of $\frac{1}{7}(1,2)$ contains two points with weights $2$ and $4$ respectively, possibly reordering indices, we may assume that $c_i=\frac{4}{7}$ when $i$ is odd and $c_i=\frac{2}{7}$ when $i$ is even. In this case, $G_i^2=-4$ when $i$ is odd and $G_i^2=-2$ when $i$ is even. Thus
$$K_Y\cdot\sum_{i=1}^sb_iF_i=\sum_{i=1}^sb_i(K_Y\cdot F_i)=\sum_{i=1}^sb_i(-2-F_i^2)\leq\sum_{i=1}^s\frac{1}{2}\cdot 2=s$$
and
\begin{align*}
        K_Y\cdot\sum_{i=1}^tc_iG_i&=\sum_{i=1}^tc_i(K_Y\cdot G_i)=\sum_{i=1}^tc_i(-2-G_i^2)=\sum_{i=1}^{\frac{t}{2}}(c_{2i-1}(-2-G_{2i-1}^2)+c_{2i}(-2-G_{2i}^2))\\
    &=\sum_{i=1}^{\frac{t}{2}}(\frac{4}{7}\cdot 2+\frac{2}{7}
    \cdot 0)=\frac{4}{7}t.\\
\end{align*}
Since $f$ extracts $k+s+t$ divisors, we have
$\rho(Y)\geq 1+k+s+t.$
Since $K_X$ is big and nef, we have $K_X^2>0$, which implies that
$$K_Y^2=K_X^2-K_Y(\cdot\sum_{i=1}^k\frac{i}{2k+1}E_i+\sum_{i=1}^sb_iF_i+\sum_{i=1}^tc_iG_i)>-\frac{k}{2k+1}-s-\frac{4}{7}t>-\frac{1}{2}-s-\frac{4}{7}t.$$
Since $X$ is rational, $Y$ is rational. By Lemma \ref{lem: k=10-rho},
$K_Y^2=10-\rho(Y).$
Thus $$-\frac{1}{2}-s-\frac{4}{7}t<K_Y^2=10-\rho(Y)\leq 10-(1+k+st)=9-k-s-t,$$
which implies that $k<\frac{19}{2}-\frac{3}{7}t<10$, a contradiction.
\end{proof}

\begin{lem}\label{lem: kx cdot c lower bound}
Then there exists a positive integer $n_1$ and a DCC set $\Ii$ of non-negative real numbers satisfying the following. Assume that
\begin{itemize}
    \item $X$ is a $\frac{2}{5}$-klt surface such that $K_X$ is big and nef, 
    \item $C$ is an irreducible curve on $X$, 
    \item $f: Y\rightarrow X$ is the minimal resolution of $X$, 
    \item $C_Y$ is the strict transform of $C$ on $Y$, and
    \item $K_Y\cdot C_Y<0$,
\end{itemize}
then 
\begin{enumerate}
    \item $K_X\cdot C\in\Ii$, and
    \item if $K_X\cdot C=0$, then $n_1K_X$ is Cartier near $C$.
\end{enumerate}
In particular, we may define $\gamma_0:=\min\{1,\gamma\in\Ii\mid \gamma>0\}.$
\end{lem}
\begin{proof}
By Lemma \ref{lem: 1/3 lc must 2 chain and 3}, there exists a positive integer $n_0=n_0(\frac{1}{15})$, such that for any closed point $X\ni x$, either $n_0K_X$ is Cartier near $x$, or $x$ is a cyclic quotient singularity of type $\frac{1}{2k+1}(1,k)$ for some positive integer $k\geq 10$. Now we let
$$\Ii:=\{\gamma\mid \gamma\geq 0,\gamma=-1+\sum_{i=1}^m\frac{k_i}{2k_i+1}+\frac{l}{n_0},m,l,k_1,\dots,k_m\in\mathbb N\}.$$
Then $\Ii$ is a DCC set of non-negative real numbers. 

Consider the equation
$$\sum_{i=1}^m\frac{k_i}{2k_i+1}+\frac{l}{n_0}=1,$$
where $m,l,k_1,\dots,k_m\in\mathbb N$. Then there exists a finite set $\Ii_0\subset\mathbb N$ such that $k_i\in\Ii_0$ for each $i$. We define
$$n_1:=n_0\prod_{\gamma\in\Ii_0}(2\gamma+1).$$

We show that $\Ii$ and $n_1$ satisfy our requirements. For any curve $C$ as in the assumption, there exists a non-negative integer $s$, such that 
\begin{itemize}
    \item there are closed points $x_1,\dots,x_s$ on $X$, such that $x_i\in C$ and $x_i$ is a cyclic quotient singularity of type $\frac{1}{2k_i+1}(1,k_i)$ for some positive integer $k_i\geq 10$ for each $i$, and
    \item for any closed point $y\not\in\{x_1,\dots,x_s\}$, $n_0K_X$ is Cartier near $y$.
\end{itemize}
By Lemma \ref{lem: cy only intersect -3 curve} and Lemma \ref{lem: 1/2 lc no k 2k+1 k >10}, we may write
$$K_Y+\sum_{i=1}^s\sum_{j=1}^{k_i}a_{i,j}E_{i,j}+\sum_{k=1}^t\frac{c_k}{n_0}F_k=f^*K_X,$$
where
\begin{itemize}
\item $E_{i,j}$ and $F_k$ are distinct prime $f$-exceptional divisors for every $i,j,k$,
\item for any $i,j$, $\Center_XE_{i,j}=x_i$,
\item $k_i,c_k$ are positive integers,
\item $a_{i,k_i}=\frac{k_i}{2k_i+1}$ for each $i$, and
\item $C_Y\cdot E_{i,u_i}=1$ and $C_Y\cdot E_{i,j}=0$ for every $j\not=u_i$. 
\end{itemize}
By Lemma \ref{lem: not nef curve is -1}, $K_Y\cdot C_Y=-1$. Thus
$$f^*K_X\cdot C_Y=(K_Y+\sum_{i=1}^s\sum_{j=1}^{k_i}a_{i,j}E_{i,j}+\sum_{k=1}^t\frac{c_k}{n_0}F_k)\cdot C_Y=-1+\sum_{i=1}^s\frac{k_i}{2k_i+1}+\frac{l}{n_0}$$
for some non-negative integer $l$. Moreover, since $K_X$ is big and nef,
$$0\leq K_X\cdot C=f^*K_X\cdot C_Y.$$
Thus $f^*K_X\cdot C_Y\in\Ii$. Moreover, if $K_X\cdot C=0$, then 
$$0=-1+\sum_{i=1}^s\frac{k_i}{2k_i+1}+\frac{l}{n_0},$$
which implies that $k_i\in\Ii_0$ for each $i$. Thus $n_1K_X$ is Cartier near $C$ by construction of $n_1$.

\end{proof}

\subsection{Construction of a nef \texorpdfstring{$\Qq$}-%
-divisors}

\begin{prop}\label{prop: find nef divisor}
There exists a positive integer $m_0$ satisfying the following. Assume that
\begin{itemize}
    \item $X$ a $\frac{2}{5}$-klt surface such that $K_X$ is big and nef, 
    \item $f: Y\rightarrow X$ is the minimal resolution of $X$, and
    \item $K_Y+\sum_{i=1}^sa_iE_i=f^*K_X$, where $E_i$ are the prime $f$-exceptional divisors,
\end{itemize}
then $m_0K_Y+\sum_{i=1}^sc_iE_i$ is nef for some non-negative integers $c_1,\dots,c_s$, such that $c_i\leq\lfloor m_0a_i\rfloor$ for each $i$.
\end{prop}
\begin{proof}
Let $n_1$ and $\gamma_0$ be the numbers given by Lemma \ref{lem: kx cdot c lower bound}, $n_2:=\max\{10,n_1,\lceil\frac{1}{\gamma_0}\rceil\},$ and 
$$m_0:=n_1\prod_{i=1}^{n_2}(2i+1).$$
We show that $m_0$ satisfies our requirements. 

We classify the singularities on $X$ into three classes:

\medskip

\noindent\textbf{Class 1}. Cyclic quotient singularities of type $\frac{1}{2k+1}(1,k)$ where $k\geq n_2$. Let these singularities be $x_1,\dots,x_s$ for some non-negative integer $s$. We may assume that $x_i$ is a cyclic quotient singularity of type $\frac{1}{2k_i+1}(1,k_i)$ for some integer $k_i\geq n_2$ for every $1\leq i\leq s$.

\medskip

\noindent\textbf{Class 2}. Singularities of type $\frac{1}{2k+1}(1,k)$ where $5\leq k<n_2$. Let these singularities be $x_{s+1},\dots,x_t$ for some integer $t\geq s$. In particular, by the definition of $m_0$. $m_0K_X$ is Cartier near $x_i$ for every $s+1\leq i\leq t$.

\medskip

\noindent\textbf{Class 3}. Other singularities.  Let these singularities be $x_{t+1},\dots,x_r$ for some integer $r\geq t$. In particular, by Lemma \ref{lem: 1/3 lc must 2 chain and 3} and the definition of $m_0$, $m_0K_X$ is Cartier near $x_i$ for every $t+1\leq i\leq r$.

\medskip

Now we may write
$$K_Y+\sum_{i=1}^s\sum_{j=1}^{k_i}\frac{j}{2k_i+1}E_{i,j}+\frac{1}{m_0}F=f^*K_X,$$
where 
\begin{itemize}
\item for every $1\leq i\leq s$ and $1\leq j\leq k_i$, $\Center_XE_{i,j}=x_i$,
\item for every $1\leq i\leq s$ and $1\leq j\leq k_i-1$, $E_{i,j}^2=-2$,
\item for every $1\leq i\leq s$, $E_{i,k_i}^2=-3$, 
\item $F\geq 0$ is a $f$-exceptional Weil divisor, such that $x_i\not\in\Center_XF$ for every $1\leq i\leq s$.
\end{itemize}
We show that we may take
$$\sum_{i=1}^lc_iE_i:=\sum_{i=1}^s\sum_{j=k_i-n_2+1}^{k_i}\frac{m_0(j-(k_i-n_2))}{2n_2+1}E_{i,j}+F.$$
Indeed, by our constructions, $0\leq c_i\leq\lfloor m_0a_i\rfloor$ for each $i$, and we only left to check that $(m_0K_Y+\sum_{i=1}^lc_iE_i)\cdot C_Y\geq 0$ for any irreducible curve $C_Y$ on $Y$. We have the following cases:

\medskip

\noindent\textbf{Case 1}. $K_Y$ is not pseudo-effective. In this case, by Lemma \ref{lem: not rational is pe}, $X$ is rational. By Lemma \ref{lem: 1/2 lc no k 2k+1 k >10}, $s=0$. Thus $\sum_{i=1}^lc_iE_i=F$ and $$(m_0K_Y+\sum_{i=1}^lc_iE_i)=m_0f^*K_X$$
is nef. Thus $(m_0K_Y+\sum_{i=1}^lc_iE_i)\cdot C_Y\geq 0$ for any irreducible curve $C_Y$ on $Y$. 

\medskip

\noindent\textbf{Case 2}. $K_Y$ is pseudo-effective.

\medskip

\noindent\textbf{Case 2.1}. $C_Y$ is not exceptional over $X$. Let $C:=f_*C_Y$.

\medskip

\noindent\textbf{Case 2.1.1}. $K_Y\cdot C_Y\geq 0$. In this case, $E_{i,j}\cdot C\geq 0$ and $F\geq C\geq 0$, and $(m_0K_Y+\sum_{i=1}^lc_iE_i)\cdot C_Y\geq 0$.

\medskip

\noindent\textbf{Case 2.1.2}. $K_Y\cdot C_Y<0$. By Lemma \ref{lem: not nef curve is -1}, $K_Y\cdot C_Y=C_Y^2=-1$. By Lemma \ref{lem: cy only intersect -3 curve}, $C_Y\cdot E_{i,j}=0$ for every $i$ and every $j\leq k_i-1$, and $C_Y\cdot E_{i,k_i}\in\{0,1\}$ for every $i$. By Lemma \ref{lem: kx cdot c lower bound}, there are two possibilities.

\medskip

\noindent\textbf{Case 2.1.2.1}. $n_1K_X$ is Cartier near $C$. In this case, since $n_2\geq n_1$, $2k_i+1\geq 2n_2+1>n_1$ for every $i$. Thus $C_Y$ does not intersect $E_{i,j}$ for any $i,j$, and hence
\begin{align*}
(m_0K_Y+\sum_{i=1}^lc_iE_i)\cdot C_Y&=(m_0K_Y+F)\cdot C_Y\\
&=(m_0K_Y+\sum_{i=1}^s\sum_{j=1}^{k_i}\frac{m_0j}{2k_i+1}E_{i,j}+F)\cdot C_Y=m_0f^*K_X\cdot C_Y\geq 0.
\end{align*}

\medskip

\noindent\textbf{Case 2.1.2.2}. $K_X\cdot C\geq\gamma_0$. Possibly reordering indices, we may assume that there exists an integer $t\in\{0,1,2,\dots,s\}$, such that $C_Y\cdot E_{i,k_i}=1$ when $1\leq i\leq t$ and $C_Y\cdot E_{i,k_i}=0$ when $t+1\leq i\leq s$. There are two cases:

\medskip

\noindent\textbf{Case 2.1.2.2.1}. $t\leq 2$. In this case, since $n_2\geq\frac{1}{\gamma_0}$, $\gamma_0>\frac{1}{2n_2+1}$. Thus
\begin{align*}
    (m_0K_Y+\sum_{i=1}^lc_iE_i)\cdot C_Y&=m_0f^*K_X\cdot C-\sum_{i=1}^t(\frac{m_0k_i}{2k_i+1}-\frac{m_0n_2}{2n_2+1})\\
    &\geq m_0\gamma_0-m_0\sum_{i=1}^t(\frac{1}{2}-\frac{n_2}{2n_2+1})\geq m_0\gamma_0-\frac{m_0}{2n_2+1}>0.
\end{align*}

\medskip

\noindent\textbf{Case 2.1.2.2.2}. $t\geq 3$. In this case, we have
\begin{align*}
    (m_0K_Y+\sum_{i=1}^lc_iE_i)\cdot C_Y&\geq m_0K_Y\cdot C_Y+\sum_{i=1}^t\frac{m_0n_2}{2n_2+1}=m_0(-1+\sum_{i=1}^t\frac{n_2}{2n_2+1})\\
    &\geq m_0(-1+\frac{3n_2}{2n_2+1})=\frac{m_0(n_2-1)}{2n_2+1}>0.
\end{align*}

\medskip

\noindent\textbf{Case 2.2}. $C_Y$ is exceptional over $X$. Then $C\subset\Supp(\cup_{i=1}^s\cup_{j=1}^{k_i}E_{i,j})\cup\Supp F.$

\medskip

\noindent\textbf{Case 2.2.1} $C_Y\subset\Supp F$. In this case, $C_Y\cdot E_{i,j}=0$ for every $i,j$, and hence
\begin{align*}
(m_0K_Y+\sum_{i=1}^lc_iE_i)\cdot C_Y&=(m_0K_Y+F)\cdot C_Y=(m_0K_Y+\sum_{i=1}^s\sum_{j=1}^{k_i}\frac{m_0j}{2k_i+1}E_{i,j}+F)\cdot C_Y\\
&=m_0f^*K_X\cdot C_Y=0.
\end{align*}

\medskip

\noindent\textbf{Case 2.2.2} $C_Y\subset\Supp(\cup_{i=1}^s\cup_{j=1}^{k_i}E_{i,j})$. We may assume that $C_Y=E_{i,j}$ for some $i$ and some $1\leq j\leq k_i$. In this case,
$$(m_0K_Y+\sum_{i=1}^lc_iE_i)\cdot C_Y=(m_0K_Y+\sum_{j=k_i-n_2+1}^{k_i}\frac{m_0(j-(k_i-n_2))}{2n_2+1}E_{i,j})\cdot C_Y.$$
There are four possibilities:

\medskip

\noindent\textbf{Case 2.2.2.1} $j=k_i$. In this case, 
$(m_0K_Y+\sum_{i=1}^lc_iE_i)\cdot C_Y=m_0(1+\frac{n_2-1}{2n_2+1}-\frac{3n_2}{2n_2+1})=0.$

\medskip

\noindent\textbf{Case 2.2.2.2} $k_i-n_2+1\leq j\leq k_i-1$. In this case, 
$$(m_0K_Y+\sum_{i=1}^lc_iE_i)\cdot C_Y=m_0(0+\frac{j-1-(k_i-n_2)}{2n_2+1}-\frac{2(j-(k_i-n_2))}{2n_2+1}+\frac{j+1-(k_i-n_2)}{2n_2+1})=0.$$

\medskip

\noindent\textbf{Case 2.2.2.3} $j=k_i-n_2$. In this case, 
$(m_0K_Y+\sum_{i=1}^lc_iE_i)\cdot C_Y=\frac{m_0}{2n_2+1}>0.$

\medskip

\noindent\textbf{Case 2.2.2.4} $1\leq j\leq k_i-n_2$. In this case, 
$(m_0K_Y+\sum_{i=1}^lc_iE_i)\cdot C_Y=m_0K_Y\cdot C_Y=0.$

\end{proof}

\begin{prop}\label{prop: find big and nef divisor}
There exists a uniform positive integer $m_2$ satisfying the following. Assume that
\begin{enumerate}
    \item $X$ a $\frac{2}{5}$-klt surface such that $K_X$ is big and nef,
    \item $f: Y\rightarrow X$ is the minimal resolution of $X$, and
    \item $K_Y+\sum_{i=1}^sa_iE_i=f^*K_X$, where $E_i$ are the prime $f$-exceptional divisors,
\end{enumerate}
then $m_2K_Y+\sum_{i=1}^lr_iE_i$ is big and nef for some non-negative integers $r_1,\dots,r_l$, such that $r_i\leq\lfloor m_2a_i\rfloor$ for each $i$.
\end{prop}

\begin{proof}
By Proposition \ref{prop: find nef divisor}, there exist a positive integer $m_0=m_0$ which does not depend on $X$, and non-negative integers $c_1,\dots,c_l$, such that $m_0K_Y+\sum_{i=1}^lc_iE_i$ is nef and $c_i\leq\lfloor m_0a_i\rfloor$ for each $i$. By Theorem \ref{thm: eb hmx14 surface}, there exists a uniform positive integer $m_1$ such that $|m_1K_X|$ defines a birational map. Let $m_2:=m_2m_1$. Then $|m_2K_X|$ defines a birational map, and hence
$$|m_2K_Y+\sum_{i=1}^l\lfloor m_2a_i\rfloor E_i|=|m_2K_Y+\sum_{i=1}^lm_2a_iE_i|=|f^*(m_2K_X)|$$
defines a birational map. 

Let $D:=m_2K_Y+\sum_{i=1}^l\lfloor m_2a_i\rfloor E_i$ and $\tilde D:=m_2K_Y+\sum_{i=1}^lm_1c_iE_i$. Since $c_i\leq\lfloor m_0a_i\rfloor$, $$m_1c_i\leq m_1\lfloor m_0a_i\rfloor\leq\lfloor m_1m_0a_i\rfloor=\lfloor m_2a_i\rfloor.$$
Thus $D\geq\tilde D$. By Proposition \ref{prop: zariski decomposition give big and nef birational map}, there exists a Weil divisor $D'$ on $X$, such that $D\geq D'\geq\tilde D$ and $D'$ is big and nef. In particular, we may write
$D'=m_2K_Y+\sum_{i=1}^lr_iE_i$
for some integers $r_1,\dots,r_l$ such that $0\leq c_i\leq r_i\leq\lfloor m_2a_i\rfloor$ for each $i$. $m_2$ and $r_1,\dots,r_l$ satisfy our requirements.
\end{proof}

\subsection{Proof of the main theorem}

\begin{proof}[Proof of Theorem \ref{thm: 1/2 lc no fixed part}]
Let $f: Y\rightarrow X$ be the minimal resolution of $X$ such that
$K_Y+\sum_{i=1}^na_iE_i=f^*K_X,$
where $E_1,\dots,E_n$ are the prime exceptional divisors of $f$. By Proposition \ref{prop: find big and nef divisor}, there exists a uniform positive integer $m_2$, such that
$K_Y+\sum_{i=1}^n\frac{r_i}{m_2}E_i$
is big and nef for some integers $r_1,\dots,r_n$ such that $0\leq r_i\leq\lfloor m_2a_i\rfloor$ for each $i$. By Lemma \ref{lem: 192 bpf}, $|192m_2K_X|$ defines a birational map and we may let $m:=192m_2$.
\end{proof}

\section{Examples}

\begin{lem}\label{lem: find smaller big and nef divisor}
Let $X$ be an lc projective surface such that $K_X$ is big and nef, $f: Y\rightarrow X$ the minimal resolution of $X$, and $E_1,\dots,E_n$ the prime $f$-exceptional divisors of $X$. Assume that
$K_Y+\sum_{i=1}^na_iE_i=f^*K_X,$
where $a_i:=1-a(E_i,X,0)$. Let $m$ be a positive integer and $c_1,\dots,c_n$ non-negative integers, such that
\begin{itemize}
\item $0\leq c_1,\dots,c_n\leq\lfloor ma_i\rfloor$,
\item $|mK_Y+\sum_{i=1}^nc_iE_i|\not=\emptyset$,
\item the fixed part of $|mK_Y+\sum_{i=1}^nc_iE_i|$ is supported on $\cup_{i=1}^nE_i$, and
\item $mK_Y+\sum_{i=1}^nc_iE_i$ is big but not nef,
\end{itemize}
then there exist non-negative integers $c'_1,\dots,c_n'$, such that
\begin{enumerate}
    \item $0\leq c_i'\leq c_i$ for each $i$,
    \item there exists $j\in\{1,2,\dots,n\}$ such that $c_j'\leq c_j$, 
    \item $|mK_Y+\sum_{i=1}^nc_iE_i|\not=\emptyset$, and
    \item the fixed part of $|mK_Y+\sum_{i=1}^nc'_iE_i|$ is supported on $\cup_{i=1}^nE_i$,
\end{enumerate}
\end{lem}
\begin{proof}
Since $0\leq c_1,\dots,c_n\leq\lfloor ma_i\rfloor$, $(Y,\sum_{i=1}^n\frac{c_i}{m}E_i)$ is lc. Thus we may run a $(K_Y+\sum_{i=1}^n\frac{c_i}{m}E_i)$-MMP $h: Y\rightarrow W$. Since the fixed part of $|mK_Y+\sum_{i=1}^nc_iE_i|$ is supported on $\cup_{i=1}^nE_i$, $h$ only contracts divisors supported on $\cup_{i=1}^nE_i$. Let $B:=h_*(K_Y+\sum_{i=1}^n\frac{c_i}{m}E_i)$, then we have
$$K_Y+\sum_{i=1}^n\frac{c_i}{m}E_i=h^*(K_W+B)+\sum_{i=1}^nb_iE_i$$
where $b_i\geq 0$ are real numbers. Moreover, since $mK_Y+\sum_{i=1}^nc_iE_i$ is big but not nef, $h\not=\id_Y$. Thus there exists $j\in\{1,2,\dots,n\}$ such that $b_j>0$. We have
$$mh^*(K_W+B)=mK_Y+\sum_{i=1}^n(c_i-mb_i)E_i.$$
Since $f$ is the minimal resolution of $X$, $E_i^2\leq -2$ for every $i$. Thus $h$ is the minimal resolution of $W$, which implies that $c_i-mb_i\geq 0$ for every $i$. Let $c_i':=\lfloor c_i-mb_i\rfloor$ for every $i$. Then (1)(2) hold. Since
$$|mK_Y+\sum_{i=1}^nc_iE_i|\cong |m(K_W+B)|\cong |mh^*(K_W+B)|\cong |mK_Y+\sum_{i=1}^nc_i'E_i|,$$
(3)(4) hold.
\end{proof}

\begin{thm}\label{thm: no fixed part implies good nef divisor}
Let $X$ be an lc projective surface such that $K_X$ is big and nef, $f: Y\rightarrow X$ the minimal resolution of $X$, and $E_1,\dots,E_n$ the prime $f$-exceptional divisors of $X$. Assume that
$K_Y+\sum_{i=1}^na_iE_i=f^*K_X$
where $a_i:=1-a(E_i,X,0)$. Then for any positive integer $m$, if $|mK_X|$ defines a birational map and does not have fixed part, then there exist positive integers $r_1,\dots,r_n$, such that
\begin{enumerate}
    \item $0\leq r_i\leq \lfloor ma_i\rfloor$, and
    \item $K_Y+\sum_{i=1}^n\frac{r_i}{m}E_i$ is big and nef.
\end{enumerate}
\end{thm}

\begin{proof}
Then the fixed part of $$|f^*(mK_X)|=|mK_Y+\sum_{i=1}^nma_iE_i|=|mK_Y+\sum_{i=1}^n\lfloor ma_i\rfloor E_i|$$
is supported on $\cup_{i=1}^nE_i$. Since $|mK_X|$ defines a birational map, $|mK_Y+\sum_{i=1}^n\lfloor ma_i\rfloor E_i|$ defines a birational map. In particular, $mK_Y+\sum_{i=1}^n\lfloor ma_i\rfloor E_i$ is big.  

We inductively define integers $c^j_i$ for every $i\in\{1,2,\dots,n\}$ for non-negative integers $j$ in the following way: Let $c^0_i:=\lfloor ma_i\rfloor$ for every $i$. If $K_Y+\sum_{i=1}^n\frac{c^j_i}{m}E_i$ is big and nef, then we let $r_i:=c^j_i$ for every $i$ and we are done. Otherwise, by Lemma \ref{lem: find smaller big and nef divisor}, there exist integers $c^{j+1}_i$ for every $i$, such that $0\leq c^{j+1}_i\leq c^j_i$,  $c^{j+1}_k<c^j_k$ for some $k\in\{1,2,\dots,n\}$, $|mK_Y+\sum_{i=1}^n\lfloor c^{j+1}_i\rfloor E_i|\not=\emptyset$, and the fixed part of $|mK_Y+\sum_{i=1}^n\lfloor c^{j+1}_i\rfloor E_i|$
is supported on $\cup_{i=1}^nE_i$. This process must terminates after finitely many steps, and we get the desired $r_i$ for every $i$.
\end{proof}

\begin{exthm}[$=$ Example-Theorem \ref{exthm: has fixed part when mld goes to 1/n}]\label{exthm: precise constuction has fixed part}
There are normal projective surfaces $\{X_{n,k}\}_{n\geq 4,k\geq 2}$, such that
\begin{enumerate}
    \item $|mK_{X_{n,k}}|\not=\emptyset$ and has non-zero fixed part for any positive integers $m,n$, and $k\geq m$.
    \item $K_{X_{n,k}}$ is ample for every $n,k$,
    \item $\lim_{k\rightarrow+\infty}\mld(X_{n,k})=\frac{1}{n-1}$
    for any $n$, and
\end{enumerate}

\begin{proof}

\noindent\textbf{Step 1}. In this step we construct $X_{n,k}$ for every $n\geq 4$ and $k\geq 2$.

For any positive integer $n\geq 4$ and positive integer $k\geq 2$, we let $Y_{n,k}$ to be a general hypersurface of degree $d_{n,k}:=2k(n-2)^2(2k(n-1)-1)$ in the weighted projective space $P_{n,k}:=\mathbb P(1,1,2k(n-2),2k(n-2)(n-1)+1)$. Since 
$2k(n-2)\mid d_{n,k}$
and
$$d_{n,k}-1=(2k(n-2)(n-1)+1)\cdot (2k(n-2)-1),$$
$Y_{n,k}$ is well-formed and has a unique singularity $o_{n,k}$, which is a cyclic quotient singularity of type 
$\frac{1}{2k(n-2)(n-1)+1}(1,2k(n-2)).$ The dual graph of this cyclic quotient singularity is the following:
	\begin{center}
		\begin{tikzpicture}[scale=1, baseline=(current bounding box.center)]
		\node[circle,draw,inner sep=0pt,minimum size=6pt,label={$2$}] (a) at (0,0){};
		\node[circle,draw,inner sep=0pt,minimum size=6pt,label={$2$}] (b) at (1.5,0){};
		\node[circle,draw,inner sep=0pt,minimum size=6pt,label={$2$},label={[shift={(-0.4,-0.35)}]$\cdots$}] (c) at (2.5,0){};
		\node[circle,draw,inner sep=0pt,minimum size=6pt,label={$2$}] (d) at (4,0){};
		\node[circle,draw,inner sep=0pt,minimum size=6pt,label={$n$}] (e) at (5.5,0){};
		
		\draw (a) edge (b);
		\draw (c) edge (d);
		\draw (d) edge (e);
		\end{tikzpicture}
	\end{center}
	where there are $2k(n-2)-1$ ``$2$" in the chain. Let $E_1=E_1(n,k),\dots,E_{2k(n-2)}=E_{2k(n-2)}(n,k)$ be the curves in this dual graph in order, i.e.
	\begin{itemize}
	    \item $E_i^2=-2$ when $i\in\{1,2,\dots,2k(n-2)-1\}$, 
	    \item $E_{2k(n-2)}^2=-n$, and
	    \item $E_i\cdot E_j\not=0$ if and only if $|i-j|\leq 1$.
	\end{itemize}
	Let $h_{n,k}:Z_{n,k}\rightarrow Y_{n,k}$ be the minimal resolution, then we have
	$$K_{Z_{n,k}}+\sum_{i=1}^{2k(n-2)}\frac{i(n-2)}{2k(n-1)(n-2)+1}E_i=h_{n,k}^*K_{Y_{n,k}}.$$
	Now let $g_{n,k}:W_{n,k}\rightarrow Z_{n,k}$ be the blow-up of $E_{k(n-1)}\cap E_{k(n-1)+1}$ and $C_{n,k,W}$ the exceptional divisor of $g_{n,k}$. Let $E_{i,W}=E_{i,W}(n,k)$ be the strict transform of $E_i$ on $W_{n,k}$ for each $i$. Then
	\begin{align*}
	    &K_{W_{n,k}}+\sum_{i=1}^{2k(n-2)}\frac{i(n-2)}{2k(n-1)(n-2)+1}E_{i,W}+\frac{n-3}{2k(n-1)(n-2)+1}C_{n,k,W}\\
	    =&g_{n,k}^*(K_{Z_{n,k}}+\sum_{i=1}^{2k(n-2)}\frac{i(n-2)}{2k(n-1)(n-2)+1}E_i)=(h_{n,k}\circ g_{n,k})^*K_{Y_{n,k}}.
	\end{align*}
	Now we run a $(K_{W_{n,k}}+\sum_{i=1}^{2k(n-2)}E_{i,W}+\frac{n-3}{2k(n-1)(n-2)+1}C_{n,k,W})$-MMP over $Y_{n,k}$ which induces a birational contraction $f_{n,k}:W_{n,k}\rightarrow X_{n,k}$. Then $f_{n,k}$ contracts precisely $E_{1,W},\dots,E_{2k(n-2),W}$. We let $C_{n,k}$ be the pushforward of $C_{n,k,W}$ on $X_{n,k}$ and $p_{n,k}: X_{n,k}\rightarrow Y_{n,k}$ the induced contraction. 
	
	\medskip
	
	\noindent\textbf{Step 2}. In this step, we show the following:
	
	\begin{claim}\label{claim: mkxnk has no fixed part when k>m}
	For any positive integers $m$, $n\geq 4$ and $k\geq m$, if $|mK_{X_{n,k}}|\not=\emptyset$ and $K_{X_{n,k}}$ is big, then $|mK_{X_{n,k}}|$ has non-zero fixed part.
	\end{claim}
	\begin{proof}
		We let $o_1=o_1(n,k):=(f_{n,k})_*(\cup_{i=1}^{k(n-1)}E_{i,W})$ and 
	and 
	$o_2=o_2(n,k):=(f_{n,k})_*(\cup_{i=k(n-1)+1}^{2k(n-2)}E_{i,W}).$
	Then $o_1$ is a cyclic quotient singularity of type
	    $\frac{1}{2k(n-1)+1}(1,k(n-1))$
	    with dual graph
	    	\begin{center}
		\begin{tikzpicture}[scale=1, baseline=(current bounding box.center)]
		\node[circle,draw,inner sep=0pt,minimum size=6pt,label={$2$}] (a) at (0,0){};
		\node[circle,draw,inner sep=0pt,minimum size=6pt,label={$2$}] (b) at (1.5,0){};
		\node[circle,draw,inner sep=0pt,minimum size=6pt,label={$2$},label={[shift={(-0.4,-0.35)}]$\cdots$}] (c) at (2.5,0){};
		\node[circle,draw,inner sep=0pt,minimum size=6pt,label={$2$}] (d) at (4,0){};
		\node[circle,draw,inner sep=0pt,minimum size=6pt,label={$3$}] (e) at (5.5,0){};
		
		\draw (a) edge (b);
		\draw (c) edge (d);
		\draw (d) edge (e);
		\end{tikzpicture}
	\end{center}
	where there are $k(n-1)-1$ ``$2$" in the chain, and $o_2$ is a cyclic quotient singularity of type
	    	    $\frac{1}{(n-3)(2k(n-1)-1)}(1,2k(n-3)-1)$
	    with dual graph
	    	\begin{center}
		\begin{tikzpicture}[scale=1, baseline=(current bounding box.center)]
		\node[circle,draw,inner sep=0pt,minimum size=6pt,label={$3$}] (a) at (0,0){};
		\node[circle,draw,inner sep=0pt,minimum size=6pt,label={$2$}] (b) at (1.5,0){};
		\node[circle,draw,inner sep=0pt,minimum size=6pt,label={$2$},label={[shift={(-0.4,-0.35)}]$\cdots$}] (c) at (2.5,0){};
		\node[circle,draw,inner sep=0pt,minimum size=6pt,label={$2$}] (d) at (4,0){};
		\node[circle,draw,inner sep=0pt,minimum size=6pt,label={$n$}] (e) at (5.5,0){};
		
		\draw (a) edge (b);
		\draw (c) edge (d);
		\draw (d) edge (e);
		\end{tikzpicture}
	\end{center}
		where there are $k(n-3)-2$ ``$2$" in the chain. Then
		$$K_{X_{n,k}}+\frac{n-3}{2k(n-1)(n-2)+1}C_{n,k}=p_{n,k}^*K_{Y_{n,k}},$$
		
		\begin{align*}
	    &K_{W_{n,k}}+\sum_{i=1}^{2k(n-2)}\frac{i(n-2)}{2k(n-1)(n-2)+1}E_{i,W}+\frac{n-3}{2k(n-1)(n-2)+1}C_{n,k,W}\\
	    =&f_{n,k}^*(K_{X_{n,k}}+\frac{n-3}{2k(n-1)(n-2)+1}C_{n,k}),
	    \end{align*}
	    and
	    	$$K_{W_{n,k}}+\sum_{i=1}^{k(n-1)}\frac{i}{2k(n-1)+1}E_{i,W}+\sum_{i=k(n-1)+1}^{2k(n-2)}\frac{i-1}{2k(n-1)-1}E_{i,W}=f_{n,k}^*K_{X_{n,k}}.$$
		    We have

	        $$f_{n,k}^*K_{X_{n,k}}\cdot C_{n,k,W}=-1+\frac{k(n-1)}{2k(n-1)+1}+\frac{k(n-1)}{2k(n-1)-1}=\frac{1}{4k^2(n-1)^2-1}<\frac{1}{35k^2}<\frac{5}{12k}.$$
	    Now for any positive even number $m=2l$, any $n\geq 4$ and any $k\geq l$, we have
	    $$\frac{\{m\cdot\frac{k(n-1)}{2k(n-1)+1}\}}{m}=\frac{\{\frac{2lk(n-1)}{2k(n-1)+1}\}}{2l}=\frac{\{-\frac{l}{2k(n-1)+1}\}}{2l}\geq\frac{5}{12l}\geq\frac{5}{12k}.$$
	    Thus for any positive even number $m=2l$, any $n\geq 4$ and any $k\geq l$,
	    $K_{W_{n,k}}+\sum_{i=1}^{2k(n-2)}\frac{c_i}{m}E_{i,W}$
	    is not nef for any integers $c_1,\dots,c_{2k(n-2)}$ such that
	    \begin{itemize}
	        \item $0\le c_i\le\lfloor\frac{mi}{2k(n-1)+1}\rfloor$ when $1\leq i\leq k(n-1)$, and
	        \item $0\le c_i\le\lfloor\frac{m(i-1)}{2k(n-1)-1}\rfloor$ when $k(n-1)+1\leq i\leq 2k(n-2)$.
	    \end{itemize}
 For any integer $n\geq 4$, any positive integer $m$ such that $|mK_{X_{n,k}}|\not=\emptyset$, and any integer $k\geq m$,
 \begin{itemize}
     \item if $K_{X_{n,k}}$ is not nef, then $|mK_{X_{n,k}}|$ has non-zero fixed part, and
     \item if $K_{X_{n,k}}$ is nef, then by Theorem \ref{thm: no fixed part implies good nef divisor}, $|mK_{X_{n,k}}|$ has non-zero fixed part.
 \end{itemize}
 	\end{proof}
	
	\medskip
	
	\noindent\textbf{Step 3}. In this step we show that $K_{X_{n,k}}$ is ample.

\begin{claim}\label{claim: kynk has no fixed part}
For any integers $n\geq 4$ and $k\geq 2$, $K_{Y_{n,k}}$ is ample, $|K_{Y_{n,k}}|$ defines a birational map, and $|K_{Y_{n,k}}|$ and has no fixed part. In particular, $|K_{Y_{n,k}}|$ defines a birational map.
\end{claim}
\begin{proof}
Let $d'_{n,k}:=d_{n,k}-\deg (-K_{P_{n,k}})$. Then
$$d'_{n,k}-(2k(n-2)(n-1)+1)=-4+2k(n-2)(n-1)(2k(n-2)-3)\geq116>0.$$
Thus $K_{Y_{n,k}}$ is ample and $|K_{Y_{n,k}}|$ defines a birational map. In particular, $|K_{Y_{n,k}}|\not=\emptyset$. 

Let $x,y,z,w$ be the coordinates of $P_{n,k}$ and $d'_{n,k}:=d_{n,k}-(1+1+2k(n-2)+(2k(n-2)(n-1)+1))$. Let $A:=(x^{d'_{n,k}}=0)$ and $B:=(y^{d'_{n,k}}=0)$. Then $A|_{Y_{n,k}}\in |K_{Y_{n,k}}|$ and $B|_{Y_{n,k}}\in |K_{Y_{n,k}}|$. We only need to show that $A|_{Y_{n,k}}\not=B|_{Y_{n,k}}$. This is equal to say that $Y_{n,k}$ does not contain the line $(x=y=0)$ in $P_{n,k}$. Suppose that $Y_{n,k}$ is defined by the homogeneous weighted polynomial $q_{n,k}(x,y,z,w)$. Since $Y_{n,k}$ is general, $z^{(n-2)(2k(n-1)-1)}\in q_{n,k}(x,y,z,w)$. Thus $Y_{n,k}$ does not contain the line $x=y=0$ and we are done.
\end{proof}

\begin{claim}\label{claim: kxnk is ample}
For any integers $n\geq 4$ and $k\geq 2$, $K_{X_{n,k}}$ is ample.
\end{claim}
\begin{proof}
For any $n,k$, by Claim \ref{claim: kynk has no fixed part}, the fixed part of $|p_{n,k}^*K_{Y_{n,k}}|$ is supported on $C_{n,k}$. Since
$$p_{n,k}^*K_{Y_{n,k}}=K_{X_{n,k}}+\frac{(n-3)}{2k(n-1)(n-2)+1}C_{n,k},$$
there exists a non-negative integer $r=r_{n,k}$ such that
$|K_{X_{n,k}}-r_{n,k}C_{n,k}|$
defines a birational map and has no fixed part. In particular, $K_{X_{n,k}}-r_{n,k}C_{n,k}$ is big and nef. If $r_{n,k}=0$, then $|K_{X_{n,k}}|\not=\emptyset$ and has no fixed part, which contradicts Claim \ref{claim: mkxnk has no fixed part when k>m}. Thus $r_{n,k}>0$.

Since $K_{X_{n,k}}+\frac{(n-3)}{2k(n-1)(n-2)+1}C_{n,k}$ is nef and big, $(K_{X_{n,k}}+\frac{(n-3)}{2k(n-1)(n-2)+1}C_{n,k})\cdot C_{n,k}=0$. Since $K_{X_{n,k}}-r_{n,k}C_{n,k}$ is nef and big, $K_{X_{n,k}}$ is nef and big and $K_{X_{n,k}}\cdot C_{n,k}>0$. In particular, $K^2_{X_{n,k}}>0$. 

For any irreducible curve $D_{n,k}$ on $X_{n,k}$ such that $D_{n,k}\not=C_{n,k}$, if $D_{n,k}\cdot C_{n,k}>0$, then
$$K_{X_{n,k}}\cdot D_{n,k}=(K_{X_{n,k}}-r_{n,k}C_{n,k})\cdot D_{n,k}+r_{n,k}C_{n,k}\cdot D_{n,k}>0,$$ 
and if $D_{n,k}\cdot C_{n,k}=0$, then
$$K_{X_{n,k}}\cdot D_{n,k}=K_{Y_{n,k}}\cdot (p_{n,k})_*D_{n,k}>0.$$
Thus $K_{X_{n,k}}$ is ample.
\end{proof}

\noindent\textbf{Step 4}. Claim \ref{claim: mkxnk has no fixed part when k>m} and Claim \ref{claim: kxnk is ample} imply (1)(2). Since $$\mld(X_{n,k})=\mld(X_{n,k}\ni o_2(n,k))=\frac{2k}{2k(n-1)-1}=\frac{1}{n-1-\frac{1}{2k}},$$
we have
$\lim_{k\rightarrow+\infty}\mld(X_{n,k})=\frac{1}{n-1}$
for any $n\geq 4$, which implies (3).
\end{proof}
\end{exthm}

\begin{exthm}\label{ex: dim 3 no fixed curve}
For any positive integer $m_0$, there exists a terminal threefold $X$ such that $K_X$ is ample but $|m_0K_X|$ is not free in codimension 2.
\end{exthm}

\begin{proof}
\noindent\textbf{Step 1}.
We start with a local construction by using the language of toric varieties. \par
Let $N=\ZZ^3$, ${\bf e}_1=(1,0,0)$, ${\bf e}_2=(0,1,0)$, ${\bf e}_3=(0,0,1)$, ${\bf w}=(1,1,0)$. Let ${\bf u}=(m,1,-b)$ and ${\bf v}=(-n,2,1)$, where $m,n,b$ are positive integers such that $nb=m+1$ and $2\nmid n$. Then all these vectors above are primitive in $N$. \par
Let $\Sigma_1$ be the fan determined by the single maximal $\text{Cone}({\bf e}_3,{\bf u},{\bf v})$. Let $X_{\Sigma_1}$ be the corresponding toric variety. Then $X_{\Sigma_1}$ is affine and the cyclic quotient singularity is of the form $\frac{1}{2m+n}(-1,2,2b+1)$. Notice that $X_{\Sigma_1}$ has an isolated singularity. \par
Let $\Sigma_2=\Sigma_1^*({\bf e}_2)$ be the star subdivision of $\Sigma_1$ at ${\bf e}_2$ (see \cite[Chapter 11]{CLS10}) and $X_{\Sigma_2}$ be the corresponding toric variety, then $g: X_{\Sigma_2}\to X_{\Sigma_1}$ is a birational morphism which is an isomorphism outside the unique toric invariant point $P\in X_{\Sigma_1}$. Since $\text{Cone}({\bf u},{\bf e}_2,{\bf v})$ is smooth, $X_{\Sigma_2}$ has only two isolated singularities, which are of type $\frac{1}{m}(1,-1,b)$ and $\frac{1}{n}(1,-1,2)$. In particular $X_{\Sigma_2}$ is terminal. We use $D_{{\bf e}_2},~D_{{\bf v}},~D_{{\bf u}},~D_{{\bf e}_3}$ to denote the corresponding toric invariant divisors. We can see that $D_{{\bf e}_2}$ is the only exceptional divisor. Let $R$ denote the proper curve in $X_{\Sigma_2}$ that corresponds to $\text{Cone}({\bf e}_2,{\bf e}_3)\in\Sigma_2$. Then $R\subset D_{{\bf e}_2}$. By \cite[Proposition 6.4.4]{CLS10}, $D_{{\bf u}}\cdot R=\frac{1}{m}$, $D_{{\bf v}}\cdot R=\frac{1}{n}$, $D_{{\bf e}_3}\cdot R=\frac{b}{m}-\frac{1}{n}>0$, and $D_{{\bf e}_2}\cdot R=-(\frac{2}{n}+\frac{1}{m})$. Therefore,
$$
0<\frac{2}{n}-\frac{b}{m}=K_{X_{\Sigma_2}}\cdot R<\frac{1}{n}
$$
\par

Let $\Sigma_3=\Sigma_2^*({\bf w})$ be the star subdivision of $\Sigma_2$ at ${\bf w}$ and $X_{\Sigma_3}$ be the corresponding toric variety, then $f: X_{\Sigma_3}\to X_{\Sigma_2}$ is a birational morphism which is an isomorphism outside the toric invariant point $Q\in X_{\Sigma_2}$ that corresponds to the maximal $\text{Cone}({\bf w},{\bf e}_2,{\bf e}_3)\in\Sigma_2$ . We use $D'_{{\bf e}_2},~D'_{{\bf v}},~D'_{{\bf u}},~D'_{{\bf e}_3},~D'_{{\bf w}}$ to denote the corresponding toric invariant divisors. Notice that $D'_{{\bf w}}$ is the only exceptional divisor of $f$ and $D'_{{\bf e}_2},~D'_{{\bf v}},~D'_{{\bf u}},~D'_{{\bf e}_3}$ are the birational transforms of $D_{{\bf e}_2},~D_{{\bf v}},~D_{{\bf u}},~D_{{\bf e}_3}$ on $X_{\Sigma_3}$. Let $R'$ denote the birational transform of $R$ on $X_3$, then $R'$ corresponds to $\text{Cone}({\bf e}_2,{\bf e}_3)\in\Sigma_3$.\par
Since ${\bf w}=\frac{1}{m}{\bf u}+\frac{b}{m}{\bf e}_3+\frac{m-1}{m}{\bf e}_2$, we have $K_{X_{\Sigma_3}}=f^*K_{X_{\Sigma_2}}+(\frac{1}{m}+\frac{b}{m}+\frac{m-1}{m}-1)D'_{{\bf w}}$, hence
$$
f^*K_{X_{\Sigma_2}}=K_{X_{\Sigma_3}}-\frac{b}{m}D'_{{\bf w}}
$$

By \cite[Lemma 6.4.2]{CLS10}, $D'_{\bf u}\cdot R'=1$. Thus for any positive integer $k$,
$$
\lf kf^*K_{X_{\Sigma_2}}\rf \cdot R'=(\frac{2}{n}-\frac{b}{m})k-\{\frac{m-kb}{m}\}
$$
\par
\noindent\textbf{Step 2}. Next we will use covering trick to make the canonical divisor ample.\par
Select a projective threefold $Z$ with the isolated quotient singularity of type $\frac{1}{2m+n}(-1,2,2b+1)$ at $P$, after resolving the singularity away from $P$ we may assume that $P$ is the only singular point on $Z$. By abuse of notations we still use $f:Y\to X$ and $g: X\to Z$ to denote the corresponding toric blow-ups defined in \textbf{Step 1}. Let $E$ be the exceptional divisor of $g$ and $R\subset E$ still be proper curve defined in \textbf{Step 1}. Then $-E$ is $g$-ample and we have $f^*K_Z-aE=K_X$, where 
$$
a=\frac{\frac{2}{n}-\frac{b}{m}}{\frac{2}{n}+\frac{1}{m}}>0
$$
\par
Let $L$ be a sufficiently ample Cartier divisor on $Z$ such that $g^*(L+K_Z)-aE$ is ample on $X$. We can find effective $A\sim 2L$ that is smooth and avoid $P$. Let $h:Z'\to Z$ be the double cover ramified along $A$. Then by the Hurwitz's Formula we have 
$K_Z'=h^*(K_Z+\frac{1}{2}A)$ and $h$ is \'etale around $P$. Let $X', Y',E',f',g'$ be the corresponding base change of $h$. Then $K_{X'}=h_X^*(K_X+\frac{1}{2}g^*A)=h_X^*(g^*(\frac{1}{2}A+K_Z)-aE)$ is ample, where $h_X: X'=X\times_Z Z'\to X$ is the canonical projection. We still defined $R$ to be the proper curve in one of the components of $E'$ as in \textbf{Step 1} and $R'$ its birational transform in $Y'$. Then we have 
$$
\lf m_0f^*K_{X'}\rf \cdot R'=(\frac{2}{n}-\frac{b}{m})m_0-\{\frac{m-m_0b}{m}\}
$$
for any positive integer $m_0$ since $Z'$ and $X_{\Sigma_1}$ are isomorphic around the isolated singular point and the computation is local. \par

Now we can select $m,n\gg 0$ such that $2bm_0<n<m$, then $(\frac{2}{n}-\frac{b}{m})m_0-\{\frac{m-m_0b}{m}\}<\frac{m_0}{n}-\frac{1}{2}<0$. Therefore $\lf m_0f^*K_{X'}\rf \cdot R'<0$, which means any effective divisor in $|m_0f^*K_{X'}|$ contains $R'$. Since $f'$ is isomorphic over the generic point of $R$, it implies that any effective divisor in $|m_0K_{X'}|$ contains $R$. 
\end{proof}

\end{document}